\documentclass[10pt]{article}
\usepackage{hyperref}
\usepackage{amssymb}
\usepackage{amsmath}

\usepackage{color}

\usepackage{latexsym}

\usepackage{amsthm}

\font\tengoth=eufm10 at 10pt
\font\sevengoth=eufm7 at 6pt
\newfam\gothfam
\textfont\gothfam=\tengoth
\scriptfont\gothfam=\sevengoth

\newcommand{\mlabel}[1]{\marginpar{#1}\label{#1}}

\newcommand{\g}{{\mathfrak g}}

\newcommand{\z}{{\mathfrak z}}

\newcommand{\fh}{{\mathfrak h}}

\newcommand{\fj}{{\mathfrak j}}

\newcommand{\fl}{{\mathfrak l}}

\newcommand{\fq}{{\mathfrak q}}
\newcommand{\fp}{{\mathfrak p}}
\newcommand{\fr}{{\mathfrak r}}

\newcommand{\ft}{{\mathfrak t}}
\newcommand{\fu}{{\mathfrak u}}

\newcommand{\fz}{{\mathfrak z}}

\renewcommand\sp{\mathfrak {sp}} 
 
\newcommand\conf{\mathfrak {conf}} 
\newcommand\heis{\mathfrak {heis}}

\newcommand{\1}{\mathbf{1}}

\newcommand{\cD}{\mathcal{D}}

\newcommand{\cH}{\mathcal{H}}

\newcommand{\cM}{\mathcal{M}}

\newcommand{\cO}{\mathcal{O}}

\newcommand{\cS}{\mathcal{S}}

\newcommand{\cW}{\mathcal{W}}

\newcommand\bx{{\bf{x}}}

\newcommand{\trile}{\trianglelefteq}
\newcommand{\subeq}{\subseteq}
\newcommand{\supeq}{\supseteq}

\newcommand{\into}{\hookrightarrow}

\newcommand{\shalf}{{\textstyle{\frac{1}{2}}}}

\newcommand{\Z}{{\mathbb Z}}
\newcommand{\R}{{\mathbb R}}
\newcommand{\C}{{\mathbb C}}

\renewcommand{\tilde}{\widetilde}

\renewcommand{\L}{\mathop{\bf L{}}\nolimits}

\newcommand{\Aff}{\mathop{{\rm Aff}}\nolimits}

\newcommand{\GL}{\mathop{{\rm GL}}\nolimits}
\newcommand{\SL}{\mathop{{\rm SL}}\nolimits}

\newcommand{\AU}{\mathop{{\rm AU}}\nolimits}

\newcommand{\PSL}{\mathop{{\rm PSL}}\nolimits}
\newcommand{\SO}{\mathop{{\rm SO}}\nolimits}

\newcommand{\OO}{\mathop{\rm O{}}\nolimits}

\newcommand{\U}{\mathop{\rm U{}}\nolimits}

\newcommand{\Sp}{\mathop{{\rm Sp}}\nolimits}
\newcommand{\Sym}{\mathop{{\rm Sym}}\nolimits}

\newcommand{\gl}  {\mathop{{\mathfrak{gl} }}\nolimits}

\newcommand{\fsl} {\mathop{{\mathfrak{sl} }}\nolimits}

\newcommand{\so}  {\mathop{{\mathfrak{so} }}\nolimits}

\newcommand{\Exp}{\mathop{{\rm Exp}}\nolimits}
\newcommand{\Fix}{\mathop{{\rm Fix}}\nolimits}

\newcommand{\ad}{\mathop{{\rm ad}}\nolimits}
\newcommand{\Ad}{\mathop{{\rm Ad}}\nolimits}

\renewcommand{\Re}{\mathop{{\rm Re}}\nolimits}
\renewcommand{\Im}{\mathop{{\rm Im}}\nolimits}

\newcommand{\tr}{\mathop{{\rm tr}}\nolimits}

\newcommand{\Pol}{\mathop{{\rm Pol}}\nolimits}

\newcommand{\Heis}{\mathop{{\rm Heis}}\nolimits}

\newcommand{\Mp}{\mathop{\rm Mp{}}\nolimits}

\newcommand{\Aut}{\mathop{{\rm Aut}}\nolimits}

\newcommand{\Conf}{\mathop{\rm Conf{}}\nolimits}

\newcommand{\Diff}{\mathop{{\rm Diff}}\nolimits}

\newcommand{\id}{\mathop{{\rm id}}\nolimits}

\newcommand{\PSO}{\mathop{{\rm PSO}}\nolimits}

\newcommand{\nin}{\noindent} 
\newcommand{\oline}{\overline}

\newcommand{\la}{\langle}
\newcommand{\ra}{\rangle}

\newcommand{\res}{\vert}

\newcommand{\Comp}{{\rm Comp}}

\newcommand{\Spec}{{\rm Spec}}

\newcommand{\ssssarr}{\hbox to 15pt{\rightarrowfill}}
\newcommand{\sssarr}{\hbox to 20pt{\rightarrowfill}}
\newcommand{\ssarr}{\hbox to 30pt{\rightarrowfill}}
\newcommand{\sarr}{\hbox to 40pt{\rightarrowfill}}
\newcommand{\arr}{\hbox to 60pt{\rightarrowfill}}
\newcommand{\sssslarr}{\hbox to 15pt{\leftarrowfill}}
\newcommand{\ssslarr}{\hbox to 20pt{\leftarrowfill}}
\newcommand{\sslarr}{\hbox to 30pt{\leftarrowfill}}
\newcommand{\slarr}{\hbox to 40pt{\leftarrowfill}}
\newcommand{\larr}{\hbox to 60pt{\leftarrowfill}}

\newcommand{\Arr}{\hbox to 80pt{\rightarrowfill}}

\def\theoremname{Theorem}
\def\propositionname{Proposition}
\def\corollaryname{Corollary}
\def\lemmaname{Lemma}
\def\remarkname{Remark}
\def\conjecturename{Conjecture} 

\def\definitionname{Definition}
\def\exercisename{Exercise}
\def\examplename{Example}
\def\examplesname{Examples}
\def\problemname{Problem}
\def\problemsname{Problems}

\def\satzname{Satz} 
\def\koroname{Korollar}
\def\folgname{Folgerung}
\def\bemerkname{Bemerkung}
\def\aufgname{Aufgabe}

\def\beisname{Beispiel}
\def\beissname{Beispiele}
\def\bewname{Beweis}

\def\@thmcounter#1{\noexpand\arabic{#1}}
\def\@thmcountersep{}
\def\@begintheorem#1#2{\it \trivlist \item[\hskip 
\labelsep{\bf #1\ #2.\quad}]}
\def\@opargbegintheorem#1#2#3{\it \trivlist
      \item[\hskip \labelsep{\bf #1\ #2.\quad{\rm #3}}]}
\makeatother
\newtheorem{theor}{\theoremname}[section]
\newtheorem{propo}[theor]{\propositionname}
\newtheorem{coro}[theor]{\corollaryname}
\newtheorem{lemm}[theor]{\lemmaname}

\newenvironment{thm}{\begin{theor}\it}{\end{theor}}

\newenvironment{prop}{\begin{propo}\it}{\end{propo}}

\newenvironment{cor}{\begin{coro}\it}{\end{coro}}

\newenvironment{lem}{\begin{lemm}\it}{\end{lemm}}

\newtheorem{rema}[theor]{\remarkname}

\newenvironment{rem}{\begin{rema}\rm}{\end{rema}}

\newtheorem{stepnow}[theor]{}

\newtheorem{defin}[theor]{\definitionname} 

\newenvironment{defn}{\begin{defin}\rm}{\end{defin}}

\newtheorem{exerc}{\exercisename}[section]

\newtheorem{exa}[theor]{\examplename}

\newenvironment{ex}{\begin{exa}\rm}{\end{exa}}

\newtheorem{exas}[theor]{\examplesname}

\newtheorem{conj}[theor]{\conjecturename}

\newtheorem{pro}[theor]{\problemname}

\newenvironment{problem}{\begin{pro}\rm}{\end{pro}}
\newenvironment{prob}{\begin{pro}\rm}{\end{pro}}

\newtheorem{prs}[theor]{\problemsname}

\newtheorem{aufg}{\aufgname}[section]

\newenvironment{prf}{\begin{proof}}{\end{proof}}
 
\newcommand{\pmat}[1]{\begin{pmatrix} #1 \end{pmatrix}}

{\hfill\qed\end{trivlist}}

\newenvironment{beweis*}{\begin{trivlist}\item[\hskip%
\labelsep{\bf\bewname.\quad}]}%
{\end{trivlist}}

\newtheorem{satzn}[theor]{\satzname}

\newtheorem{koro}[theor]{\koroname}

\newtheorem{folg}[theor]{\folgname}

\newtheorem{bem}[theor]{\bemerkname}

\newtheorem{aufgn}[theor]{\aufgname}

\newtheorem{beis}[theor]{\beisname}

\newtheorem{beiss}[theor]{\beissname}

\newcommand{\csp}{{\mathfrak{csp}}} 
\newcommand{\hcsp}{{\mathfrak{hcsp}}}

\renewcommand{\phi}{\varphi}

\newcommand{\Stand}{\mathop{{\rm Stand}}\nolimits}

\newcommand{\sV}{{\tt V}}

\renewcommand\mlabel{\label}

\begin{document}

\title{Semigroups in 3-graded Lie groups and \\ 
endomorphisms of standard subspaces}
\author{Karl-Hermann Neeb} 

\maketitle

\abstract{Let $\sV$ be a standard subspace in the complex 
Hilbert space $\cH$ and $U \colon G \to \U(\cH)$ 
be a unitary representation of a finite dimensional Lie group. 
We assume the existence of an element 
$h \in \g$ such that $U(\exp th) = \Delta_\sV^{it}$ 
is the modular group of $\sV$ and that the modular involution 
$J_\sV$ normalizes $U(G)$. We want to determine the semigroup 
$S_\sV = \{ g\in G \colon  U(g)\sV \subeq \sV\}.$ 
In previous work we have seen that its infinitesimal 
generators span a Lie algebra on which $\ad h$ defines a $3$-grading,
and here we completely determine the semigroup $S_\sV$ under the assumption 
that $\ad h$ defines a $3$-grading on $\g$. 
Concretely, we show that the $\ad h$-eigenspaces $\g^{\pm 1}$ contain 
closed convex cones $C_\pm$, such that 
\[ S_\sV = \exp(C_+) G_\sV \exp(C_-),\] 
where $G_\sV = \{ g \in G \colon U(g) \sV = \sV \}$ is the stabilizer 
of $\sV$. To obtain this result we compare several subsemigroups 
of $G$ specified by the grading and the positive cone $C_U$ of $U$. In particular, 
we show that the orbit $\cO_\sV = U(G)\sV$ with the inclusion order 
is an ordered symmetric space covering the adjoint orbit 
$\cO_h =  \Ad(G)h$, endowed with the partial order defined by~$C_U$.
}

\tableofcontents 

\section{Introduction} 
\mlabel{sec:1}

A real subspace $\sV$ of a complex Hilbert space $\cH$ is called 
{\it standard} if it is closed and 
\begin{equation}
  \label{eq:stansub}
  \sV \cap i \sV = \{0\}\quad \mbox{ and } \quad \cH = \oline{\sV + i \sV}
\end{equation} 
(cf.\ \cite{Lo08} for the basic theory of standard subspaces). 
If $\sV \subeq \cH$ is a standard subspace, then 
\begin{equation}
  \label{eq:tomitaop}
 T_\sV \colon \cD(T_\sV) := \sV + i \sV \to \cH, \quad 
x + i y \mapsto x- iy 
\end{equation}
defines a closed operator with $\sV = \Fix(T_\sV)$, called the 
\textit{Tomita operator} of~$\sV$. Its polar decomposition can be written as 
$T_\sV = J_\sV \Delta_\sV^{1/2}$, 
where $J_\sV$ is a {\it conjugation} (an antiunitary 
involution) and $\Delta_\sV$ is a positive selfadjoint operator such that the 
{\it modular relation} 
\begin{equation}
  \label{eq:modrel}
 J_\sV \Delta_\sV J_\sV = \Delta_\sV^{-1} 
\end{equation}
holds. We call  $(\Delta_\sV, J_\sV)$ the {\it pair of modular 
objects} associated to~$\sV$. 

Standard subspaces arise naturally in the modular theory 
of von Neumann algebras. 
If $\cM \subeq B(\cH)$ is a von Neumann algebra 
and $\Omega \in \cH$ is  
{\it cyclic} for $\cM$ ($\cM\Omega$ is dense in $\cH$) and 
{\it separating} (the map $\cM \to \cH,M \mapsto M\Omega$ is injective), 
then the  Tomita--Takesaki Theorem (\cite[Thm.~2.5.14]{BR87}) 
implies that $\sV := \oline{\{ M\Omega \colon M = M^* \in \cM\}}$ 
is standard, and that 
\[ J_\sV \cM J_\sV = \cM' \quad \mbox{ and } \quad 
\Delta_\sV^{it} \cM \Delta_\sV^{-it}  = \cM \quad \mbox{ for }  \quad t \in \R.\] 
So we obtain a one-parameter group of automorphisms of $\cM$ 
(the modular group) 
and a symmetry between $\cM$ and its commutant $\cM'$, implemented by $J_\sV$. 

Building on the  Haag--Kastler theory 
of local observables in Quantum Field Theory (QFT) 
(\cite{Ha96}, \cite{BS93}, \cite{BDFS00}, \cite{FR19}), 
the current interest in standard subspaces arose 
in the 1990s from the work of Borchers and Wiesbrock (\cite{Bo92, Wi93}). 
This in turn led to modular localization 
in Quantum Field Theory introduced by  
Brunetti, Guido and Longo in \cite{BGL02, BGL94, BGL93}; 
see also \cite{Le15, LL15} for important applications of this technique. 

The order on the set $\Stand(\cH)$ of standard subspaces of 
$\cH$ is of particular importance because it reflects 
inclusions of corresponding von Neumann algebras 
(see~\cite[\S 4.2]{NO17}, \cite{Lo08} and \cite{Ta10} 
for more details on the translation process). 
As the order on $\Stand(\cH)$ is hard to understand, it is 
natural to probe the ordered space $\Stand(\cH)$ by finite dimensional 
homogeneous submanifolds arising as orbits under unitary 
representations of finite dimensional Lie groups. To link such a 
representation as closely as possible to standard subspaces, 
we consider the following setting. 
Let $\sV \subeq \cH$ be a standard subspace 
and $U \colon G \to \U(\cH)$ be a unitary representation of the connected 
Lie group $G$. We further assume that there exists an involutive automorphism 
$\tau_G$ of $G$ and $h \in \g$ fixed by $\tau = \L(\tau_G)$ such that $U$ extends 
to an antiunitary representation of $G \rtimes \{\id_G,\tau_G\}$ 
such that 
\begin{equation}
  \label{eq:moddat}
  U(\exp th) = \Delta_\sV^{-it/2\pi} \quad \mbox{ for } \quad t \in \R 
\quad \mbox{ and } \quad 
J_\sV U(g) J_\sV = U(\tau_G(g)) \quad \mbox{ for } \quad g \in G.
\end{equation}
Then the order on the orbit $\cO_\sV := U(G)\sV$ is determined by 
the subsemigroup 
\begin{equation}
  \label{eq:SM=SV}
S_\sV :=  \{ g\in G \colon U(g)\sV \subeq \sV\}.
\end{equation}
In \cite{Ne19} we managed to calculate its {\it Lie wedge} 
\begin{footnote}{In the theory of 
Lie semigroups (\cite{HHL89, HN93}) Lie wedges are the semigroup analogs 
of the Lie algebras of  closed subgroups. A {\it Lie wedge} is a closed convex 
cone $W$ in a Lie algebra $\g$ such that 
$e^{\ad x} W = W$ for $x \in W \cap -W$. In particular, 
linear subspaces are Lie wedges if and only if they are Lie subalgebras.}
\end{footnote}
\[ \L(S_\sV) = \{ x \in \g \colon 
\exp(\R_+ x) \subeq S_\sV\}, \] 
i.e., the set of generators of its one-parameter 
subsemigroups in the Lie algebra $\g$ of $G$. 
To formulate this result, for $x \in \g$, we write 
$\partial U(x)$ for the (possibly unbounded) skew-adjoint operator on $\cH$ with 
$U(\exp tx) = e^{t\partial U(x)}$ for $t \in \R$, and 
\[ C_U := \{ x \in \g \colon - i \partial U(x) \geq 0\} \] 
for the {\it positive cone of $U$}.
The Structure Theorem (\cite[Thm.~4.4]{Ne19}) asserts that, 
under the natural assumption that $\ker(U)$ is discrete,  
\begin{equation}
  \label{eq:(b)}
\L(S_\sV) = C_- \oplus \g_\sV \oplus C_+, \quad \mbox{ where }  \quad 
\g_\sV = \L(G_\sV), \quad G_\sV = \{g \in G \colon U(g)\sV = \sV\}, 
\end{equation}
and 
\[ C_\pm := \pm \{ x \in C_U \colon \tau(x) = -x, [h,x] = \pm x\}. \] 
We also show that the cone $\L(S_\sV)$ spans a Lie subalgebra which is $3$-graded by $\ad h$ 
and on which $\tau$ coincides with $e^{\pi i \ad h}$. 

In the present paper we therefore focus on the situation 
where $\g$ is $3$-graded by $\ad h$ in the sense that the $\ad h$-eigenspaces 
$\g^\lambda := \g^\lambda(h) := \ker(\ad h - \lambda \1)$ satisfy 
  \begin{equation}
    \label{eq:threegrad}
\g = \g^{-1} \oplus \g^0 \oplus \g^1,  
\quad \mbox{ and  }  \quad  \tau = e^{\pi i \ad h}.  
  \end{equation}
The $\tau$-eigenspaces are then 
$\fh = \g^0$ and $\fq = \g^1 \oplus \g^{-1}$. 
We also assume that $C \subeq \g$ is a closed pointed 
$\Ad(G)$-invariant convex cone satisfying $\tau(C) = - C$ and that 
$G$ is a $1$-connected Lie group with Lie algebra~$\g$. 
The involution $\tau$ then integrates to an involution $\tau_G$ on $G$ 
and 
\[ C_\pm := \pm C \cap \g^{\pm 1} \] 
are pointed convex cones invariant under the action of the group 
$G^0 := \{ g \in G \colon \Ad(g)h = h\}$. 
These structures lead to three subsemigroups of $G$: 
\begin{itemize}
\item The Olshanski semigroup 
$S(C_\fq) := G^0 \exp(C_\fq)$ for $C_\fq := C_+ \oplus C_- \subeq \fq$, 
\item the semigroup  
$S(h,C) := \{ g \in G \colon  h - \Ad(g)h \in  C\}$, and 
\item using the complex Olshanski semigroup 
  $S(iC) := G \Exp(iC)$ (see Subsection~\ref{subsec:2.3} for details), 
we define the  
subsemigroup $G_{i\pi}(C)$ of those elements $g \in G$ 
for which the orbit map 
$\beta^g \colon \R \to G, \beta^g(t) = \exp(th) g\exp(-th)$ 
extends analytically to a map 
$\beta^g \colon \oline{\cS_{\pi}} \to S(iC)$.
\end{itemize}
The main results of Section~\ref{sec:2} are the equalities 
\begin{equation}
  \label{eq:equals}
 S(C_\fq) = S(h,C) = \exp(C_+) G^0 \exp(C_-)= G_{\pi i}(C).
\end{equation}
The first two equalities constitute the Decomposition Theorem~\ref{thm:1.7},
and the last equality is Theorem~\ref{thm:1.29}. 
The key point of the identity \eqref{eq:equals} is that it provides 
three rather different perspectives on the same subsemigroup 
of $G$, and this contains important information on the semigroups~$S_\sV$. 

To see this connection, let us first consider an 
antiunitary representation $(U,\cH)$ with discrete kernel for a 
semidirect product $G \rtimes \{\id_G, \tau_G\}$, 
where $G$ is a connected Lie group, and 
$\tau_G$ is an involutive automorphism of $G$. 
We consider the standard subspace  $\sV \subeq \cH$ specified by 
the modular objects 
\begin{equation}
  \label{eq:modobrel}
J_\sV = U(\tau_G) \quad \mbox{ and } \quad 
\Delta_\sV = e^{2\pi i \partial U(h)} \quad \mbox{ for some }\quad h \in \g^\tau 
\end{equation}
but make no further assumptions on $h$. 
Building on some observations by  Borchers and Wiesbrock 
(Theorem~\ref{thm:b.2}), we show in the 
Monotonicity Theorem (Theorem~\ref{thm:inc}) that 
\begin{equation}
  \label{eq:montheo}
S_\sV \subeq S(h,C_U).
\end{equation}
Here the main point is that, for two standard subspaces 
$\sV_1 \subeq \sV_2$, we have 
$\log \Delta_{\sV_2} \leq \log \Delta_{\sV_1}$ 
in the sense of quadratic forms. Since these logarithms 
are typically not semibounded, the order relation 
requires some explanation that we provide in Appendix~\ref{app:a}. 
Put differently, the Monotonicity Theorem asserts that the well-defined  
$G$-equivariant map 
\[ \cO_\sV = U(G)\sV \cong G/G_\sV \to \cO_h \cong G/G^0, \quad U(g)\sV \mapsto \Ad(g)h \] 
is monotone, hence the name. 

Combining the Monotonicity Theorem with the identities 
\eqref{eq:equals}, it is now easy to determine the 
semigroup $S_\sV$ for the case where $\g$ is $3$-graded by $\ad h$, and 
$\tau = e^{\pi i \ad h}$. It is given by 
\[ S_\sV = G_\sV \exp(C_\fq) = \exp(C_+) G_\sV \exp(C_-) \subeq S(h,C)
\quad \mbox{ for } \quad C = C_U.\]  
As a consequence for the general case, the infinitesimally 
generated subsemigroup $\la \exp(\L(S_\sV)\ra$ of $S_\sV$ 
coincides with $\exp(C_+) (G_\sV)_0 \exp(C_-)$. 

We conclude this paper with a short section on perspectives, 
where we explain how the present results can be 
used to explore covariant nets of standard subspaces 
on the abstract level (\cite{MN20}) and how to find them 
in Hilbert spaces of holomorphic functions on tubes 
(\cite{NOO20}) and in Hilbert spaces of distributions on Lie groups 
(\cite{NO20}). Another important issue is the classification 
of all tuples $(\g,h,C)$, where 
$\g$ is $3$-graded by $\ad h$ and generated by $C_\pm$ and $h$, 
and, more generally, the subalgebras generated by 
$C_\pm$ and $\fz_\fh(h)$ for tuples $(\g,h,\tau,C)$, where 
$\tau(h) = h$, $C$ is a pointed invariant convex cone satisfying 
$\tau(C) =- C$ and $C_\pm = \pm C \cap \fq^{\pm 1}(h)$. 
For first steps in this classification program we refer to \cite{Oeh20}. 

Finally, we note that Longo and Witten obtain in \cite{LW11} 
some results on semigroups $S_\sV$ that can be interpreted as infinite dimensional 
versions, where $G \cong E \rtimes_\alpha  \R$, and $E$ is a 
topological vector space. Their results show that, extending 
our results to infinite dimensional groups requires completely new techniques 
that are different from what we use in the finite dimensional 
case, here and in \cite{Ne19}.

\subsection*{Notation}

\begin{itemize}
\item For a Lie group $G$, we write $\g$ for its Lie algebra, 
$\Ad \colon G \to \Aut(\g)$ for the adjoint action of $G$ on $\g$, induced by the 
conjugation action of $G$ on itself, and $\ad x(y) = [x,y]$ for the adjoint 
action of $\g$ on itself. 
\item $\AU(\cH)$ is the group of unitary or antiunitary operators on a complex 
Hilbert space $\cH$. 
\item If $\tau_G \in \Aut(G)$ is an order two automorphism, 
then $G_\tau := G \rtimes \{\id_G, \tau_G\}$ becomes 
a Lie group and an {\it antiunitary representation} 
of $G_\tau$ is a homomorphism 
$U \colon G \to \AU(\cH)$ for which $J := U(\tau_G)$ 
is antiunitary and $U(G) \subeq \U(\cH)$. 
Antiunitary representations are 
assumed to be continuous with respect to the strong operator topology on 
$\AU(\cH)$. 

The automorphism of $\g$ induced by $\tau_G$ is denoted $\tau$ 
and we write 
\begin{equation}
  \label{eq:hq}
 \fh := \{ x \in \g \colon \tau(x) = x \} \quad \mbox{ and } \quad 
\fq := \{ x \in \g \colon \tau(x) = -x \} 
\end{equation}
for its eigenspaces. We put $g^\sharp := \tau_G(g)^{-1}$, 

\item For a real standard subspace $\sV \subeq \cH$, we write 
$(\Delta_\sV, J_\sV)$ for the corresponding pair 
of {\it modular objects} specified by $\sV = \Fix(J_\sV \Delta_\sV^{1/2})$. 
\item Horizontal strips in the complex plane are denoted 
$\cS_{\beta} := \{ z \in \C \colon 0  < \Im z < \beta\}$ for $\beta > 0$. 
\item For a unitary representation $U \colon G \to \U(\cH)$ of a finite 
dimensional Lie group~$G$, 
we write $\cH^\infty$ for the dense subspace of {\it smooth vectors} $\xi$, 
i.e., the orbit map $U^\xi \colon G \to \cH, g \mapsto U_g \xi$ is smooth. 
The infinitesimal generator of the unitary one-parameter group 
$(U(\exp tx))_{t \in \R}$ is denoted $\partial U(x)$. 
The closed convex $\Ad(G)$-invariant cone  
\[ C_U := \{ x \in \g \colon -i \partial U(x) \geq 0 \} \] 
is called the {\it positive cone} of the representation~$U$. 
\end{itemize}

\section{Groups and semigroups in $3$-graded Lie groups} 
\mlabel{sec:2}

In this section we take a closer look at the groups and 
semigroups that arise naturally in our setting. 
Throughout, $G$ is a $1$-connected (connected and 
simply connected) Lie group with Lie algebra $\g$,  
and $h \in \g$ defines a $3$-grading 
$\g = \g^{-1} \oplus \g^0 \oplus \g^1$ 
(see \eqref{eq:threegrad}). 
In Subsection~\ref{subsec:2.1} we study basic properties of the subgroups 
\[ G^{\pm 1} := \exp(\g^{\pm 1}), \qquad G^0 =  
 \{g \in G \colon \Ad(g)h = h\} \quad \mbox{ and }\quad 
P_\mp := G^0 G^{\mp 1}. \] 

In Subsection~\ref{subsec:2.3}, we then also 
take an $\Ad(G)$-invariant pointed closed convex cone $C$ into account 
which also satisfies $\tau(C) = -C$ for $\tau = e^{\pi i \ad h}$. 
First, we show that we have an Olshanski semigroup 
\[ S(C_\fq) := G^0 \exp(C_\fq) \quad \mbox{ for } \quad 
C_\fq := C_+ \oplus C_-.\] 

In Subsection~\ref{subsec:2.2} we use the pointed invariant cone $C$ 
to define an $\Ad(G)$-invariant
 partial order on $\g$ by $x \leq_C y$ if $y - x \in C$. 
Clearly, each adjoint orbit thus inherits an invariant 
order structure, and the orbit $\cO_h := \Ad(G)h$ 
is of particular interest. It is an ordered 
symmetric space, which leads us to the semigroup 
\[ S(h,C) := \{ g \in G \colon \Ad(g)h \leq_C h \}.\] 
One of our main results is Theorem~\ref{thm:1.7}, 
asserting that $S(C_\fq) = S(C,h)$ and that 
this semigroup has a triangular decomposition 
$S(h,C) = \exp(C_-) G^0 \exp(C_+)$ 
(Subsection~\ref{subsec:2.2}). 
In Theorem~\ref{thm:1.29} we further show that this 
subsemigroup also coincides with $G_{\pi i}(C)$. 
So we obtain three rather different perspectives on the same 
subsemigroups of~$G$ (Subsection~\ref{subsec:2.6}).

\subsection{Subgroups associated with the $3$-grading} 
\mlabel{subsec:2.1}


\begin{lem} \mlabel{lem:tridec-solrad}
If $G$ is a $1$-connected Lie group and $\g$ is $3$-graded by the 
element $h \in \g$, then there exists a 
Levi complement $\fl \subeq \g$ invariant under $\ad h$. We 
write $L$ for the corresponding integral subgroup 
and $R$ for the solvable radical of $G$. 
Then the following assertions hold: 
\begin{itemize}
\item[\rm(i)] $G \cong R \rtimes L$ 
and $G^j = R^j \rtimes L^j$ for $j = -1,0,1$. 
\item[\rm(ii)] $R = R^{-1} R^0 R^1$. 
\item[\rm(iii)] The projection $p_L \colon G \to L$ satisfies 
$G^{-1} G^0 G^1 = p_L^{-1}(L^{-1} L^0 L^1).$ 
\end{itemize}
\end{lem}

\begin{prf} Since the derivation $\ad h$ is semisimple, 
\cite[Thm.~B.2]{Ne19} implies the existence of a Levi complement $\fl$ in 
$\g$, invariant under $\ad h$. 

\nin (i)  As $G$ is simply connected, 
we have $G \cong R \rtimes L$, where $R$ and $L$ are both simply connected. 
Since both factors are 
invariant under conjugation with $\exp(th)$ for $t \in\R$, 
and $g \in G^0$ is equivalent to 
$g \exp(th) g^{-1} = \exp(th)$ for all $t\in \R,$ 
which in turn is equivalent to 
$g  = \exp(th) g \exp(th)^{-1}$ for $t\in \R,$ 
we obtain $G^0 = R^0 \rtimes L^0$. 
The relation $G^j = R^j \rtimes L^j$ for $j = \pm 1$ follows from 
$\g^{\pm 1} = \fr^{\pm 1} \rtimes \fl^{\pm 1}$, which follows from the 
invariance of $\fr$ and $\fl$ under $\ad h$. 

\nin (ii) As $\fr^{\pm 1} \subeq [h,\fr] \subeq [\g,\fr]$ 
is contained in the maximal 
nilpotent ideal~$\fu$ of $\g$ (\cite[Ch.~1, \S 5.3, Thm.~1]{Bou90}), 
we have $\fr = \fu + \fr^0(h)$. 
Hence it suffices to show that the multiplication map 
$U^{-1} \times U^0 \times U^1 \to U$ of the integral 
subgroup $U$ corresponding to $\fu$ 
is a diffeomorphism. As $\fu^{\pm 1}$ and $\fu^{\pm 1} + \fu^0$ are 
Lie subalgebras of the nilpotent Lie algebra $\fu$, this 
follows by applying \cite[Lemma~11.2.13]{HN12} twice. 

\nin (iii) As $p_\fl \colon \g \to \fl$ is a morphism of $3$-graded Lie algebras, 
we have $p_\fl(\g^j)=  \fl^j$ for $j = 0, \pm 1$, and therefore 
$p_L(G^{\pm 1}) = L^{\pm 1}$. To see that $p_L(G^0) \subeq L^0$, we note that 
$g \in G^0$ means that $\Ad(g)h = h$, and applying $p_\fl$ leads to 
\[ h_\fl := p_\fl(h) = p_\fl(\Ad(g)h) = \Ad(p_L(g)) h_\fl.\] 
If, conversely, $g_L \in L$ satisfies $\Ad(g_L) h_\fl = h_\fl$, 
then $h_\fr := h - h_\fl \in \fr$ satisfies 
$[h_\fr,\fl] = \{0\}$. As $L$ is connected, it follows that 
$\Ad(g_L)h_\fr = h_\fr$. 
We conclude that $\Ad(g_L)h = h$, and thus $p_L(G^0) = L^0$. 
This shows that 
\begin{equation}
  \label{eq:projtridec}
p_L(G^{-1}G^0 G^1) = L^{-1} L^0 L^1.
\end{equation}
Further, $R = R^{-1} R^0 R^1$ by (ii), so that the subgroup property  of 
$G^0 G^1$ and 
\eqref{eq:projtridec} lead to   
\[ p_L^{-1}(L^{-1} L^0 L^1) 
= R G^{-1} G^0 G^{1} 
=  G^{-1}  R G^0 G^{1} 
=  G^{-1}  R^{-1} R^0 R^1 G^0 G^{1} 
=  G^{-1}  G^0 G^{1}.\qedhere \] 
\end{prf}

The following lemma is useful for reductions from $G$ to its adjoint group. 

\begin{lem} \mlabel{lem:2grps}
The subgroup $G^0 = \{ g \in G \colon  \Ad(g)h = h \}$ coincides with 
\[ G_{\ad h} 
:= \{g \in G \colon \Ad(g) \ad h \Ad(g)^{-1} = \ad h\} 
= \{g \in G \colon \Ad(g) h - h \in \fz(\g)\}.\] 
\end{lem}

\begin{prf} As $G^0 \subeq G_{\ad h}$, we 
have to show that, if $\Ad(g)$ commutes with $\ad h$, i.e., if 
it preserves the $3$-grading, then $\Ad(g)h = h$. 
Let $\g = \fr \rtimes \fl$ be an $\ad h$-invariant Levi decomposition 
(Lemma~\ref{lem:tridec-solrad}) 
and write $h = h_\fr + h_\fl$, accordingly. Then $[h_\fr,\fl] = \{0\}$ 
because $\ad h$ and $\ad h_\fl$ have the same restriction on $\fl$. 
We write $G = R \rtimes L$ for the corresponding decomposition of~$G$ 
and, accordingly, $g \in G$ as $g = g_R g_L$ with $g_R \in R$, $g_L \in L$.

Assume that $\Ad(g)$ commutes with $\ad h$. 
Then $\Ad(g_L)$ commutes with $\ad h_\fl$, and this implies that 
$\Ad(g_L)h_\fl = h_\fl$ because $\fz(\fl) = \{0\}$. 
For $g = g_R g_L$ we also have 
$\Ad(g_L) h_\fr = h_\fr$ because $h_\fr$ commutes with $\fl$. 
We thus obtain 
\[ \Ad(g) h = \Ad(g_R) h \in h + \fz(\g).\] 
Next we write $g_R = \exp(x_1) \exp(x_{-1}) g_0$ with $g_0 \in R^0$ and 
$x_{\pm 1} \in \fr^{\pm 1}$ (Lemma~\ref{lem:tridec-solrad}(ii)) and obtain 
\[ \Ad(g_R) h 
= e^{\ad x_1} e^{\ad x_{-1}} h 
= e^{\ad x_1} (h + [x_{-1},h]) 
= e^{\ad x_1} (h + x_{-1}) \in h + x_{-1} + \fr^0 + \fr^1.\] 
Therefore $\Ad(g)h \in h + \fz(\g) \subeq h + \fr^0$ implies $x_{-1} = 0$. 
We likewise obtain $x_1 = 0$, so that $\Ad(g)h = \Ad(g_R) h = h$.   
\end{prf}


\begin{rem} \mlabel{rem:1.2} (a) The subgroup $G^0$ contains the center of $G$ 
but $\tau_G$ may act non-trivially on the center. A typical example 
arises for $G =\tilde\SL_2(\R)$ with 
$Z(G) \cong \Z \cong \pi_1(\PSL_2(\R))$. 
Here the fundamental group of $\PSL_2(\R)$ is generated by the loop 
obtained from the inclusion 
$\PSO_2(\R) \into \PSL_2(\R)$. 
Let $\tau_G \in \Aut(G)$ be the involution given on the Lie algebra level by 
\begin{equation}
  \label{eq:tsl2}
\tau\pmat{a & b \\ c & d} =   \pmat{a & -b \\ -c & d}, 
\quad \mbox{ which is } \quad 
e^{\pi i \ad h} \quad \mbox{ for } \quad 
h = \frac{1}{2}\pmat{1 & 0 \\ 0 & -1}.
\end{equation}
Then $\tau_G$ induces the inversion on $Z(G)$. 
Here $\exp (\R h)$ is the identity component of $G^0$, 
$Z(G) \subeq G^0$, and $Z(G) \cap \exp(\R h) = \{e\}$ so that 
$Z(G) \cong \pi_0(G^0)$. 
In particular, $G^0$ has infinitely many connected  components. 

\nin (b) Let $R \trile G$ denote the solvable radical, i.e., the maximal connected 
normal solvable subgroup. Then Levi's Theorem, and the $1$-connectedness of $G$ 
imply that $G \cong R \rtimes L$ for a $1$-connected 
semisimple Lie group~$L$. 
Then \cite[Thm.~B.2(ii)]{Ne19} implies that $R^0 = G^0 \cap R$ 
is connected, but $G^0$ need not be connected, as we have seen in (a) 
for $G = \tilde\SL_2(\R)$.

\nin (c) Let 
\[ G^\sharp := \{ g \in G \colon  \Ad(g)\tau \Ad(g)^{-1} = \tau\} 
= \Ad^{-1}(\Ad(G)^\tau)\supeq G^0.\]
 Then $g \in G^\sharp$ is equivalent to 
$g \tau(g)^{-1} = gg^\sharp \in \ker \Ad = Z(G)$. The map 
\[ \gamma \colon G^\sharp \to Z(G), \quad 
\gamma(g) = gg^\sharp \] 
is a group homomorphism because 
$\gamma(g) \gamma(h) = gg^\sharp hh^\sharp = gh h^\sharp g^\sharp  
= \gamma(gh)$ for $g, h \in G^\sharp$ follows from $hh^\sharp \in Z(G)$. 
As $G$ is $1$-connected, $G^{\tau_G} = \ker \gamma$ is connected 
(\cite[Thm.~B.2]{Ne19}), 
so that $\pi_0(G^\sharp) \cong \gamma(G^\sharp)$ can be identified with a 
subgroup of $Z(G)$.

For $g \in G^\sharp$ we have $\gamma(g)^\sharp = \gamma(g)$, i.e., 
$\tau(\gamma(g)) = \gamma(g)^{-1}$. Therefore $\gamma(g)$ is contained in the 
discrete subgroup 
\[ Z_\sharp 
:= \{ z \in Z(G) \colon z^\sharp = z\} = \{ z \in Z(G) \colon \tau(z) = z^{-1}\}.\] 
The discreteness of this subgroup follows from 
$\L(Z_\sharp) = \fz(\g) \cap \fq = \{0\}$. We also note, since 
$G$ is $1$-connected, 
$\exp\res_{\fz(\g)} \colon \fz(\g) \to Z(G)_0$ 
is bijective, so that $Z(G)_0 \cap Z_\sharp = \{e\}$. 
\end{rem}

\begin{lem} \mlabel{lem:conefield-a}
For the subgroups $P_\pm := G^0 G^{\pm 1}$, the following assertions hold: 
\begin{itemize}
\item[\rm(i)] 
$P_\pm = \{g \in G \colon \Ad(g) h \in h + \g^{\pm 1} \}$, and both 
subgroups are closed. 
\item[\rm(ii)] $P_+ \cap P_- = G^0$. 
\item[\rm(iii)] The subalgebras $\fp_\pm = \g^0 \oplus \g^{\pm 1}$ are 
self-normalizing. 
\item[\rm(iv)] The exponential function 
$\exp \colon \g^{\pm 1} \to G^{\pm 1}$ is a global diffeomorphism. 
\item[\rm(v)] The multiplication map 
$G^1 \times G^0 \times G^{-1} \to G, 
(g_1, g_0, g_{-1}) \mapsto g_1 g_0 g_{-1}$ 
is a diffeomorphism onto the  open subset $G^1 G^0 G^{-1}$ of $G$. 
\end{itemize}
\end{lem}

\begin{prf} (i) As $\Ad(P_\pm) h = \Ad(G^{\pm 1})h 
= h + \g^{\pm 1}$, we have $\Ad(g)h - h \in \g^{\pm 1}$ for $g \in P_\pm$. 
If, conversely, $g \in G$ satisfies 
$x_{\pm 1} := \Ad(g)h - h  \in \g^{\pm 1}$, 
then 
\[ \Ad( (\exp \pm x_{\pm 1})g)h 
= e^{\pm \ad x_{\pm 1}}(h + x_{\pm 1})
= h + x_{\pm 1} \pm  [x_{\pm 1}, h] 
= h + x_{\pm 1} - x_{\pm 1} = h\] 
implies that $\exp(\pm x_{\pm 1}) g \in G^0$, so that $g \in G^{\pm 1} G^0 = P_{\pm}$. 

\nin (ii) Clearly, $G^0 \subeq P_+ \cap P_-$. If, conversely, 
$g \in P_+ \cap P_-$, then 
$\Ad(g) h -h \in \g^{1} \cap \g^{-1} = \{0\}$ by (i), so that $g \in G^0$. 

\nin (iii) In view of $h \in \g^0$, for $x = x_{-1} + x_0 + x_1$ with 
$x_j \in \g^j$, the relation 
$x_1 - x_{-1} = [h,x] \in \fp_\pm$ implies $x_{\mp 1} = 0$, so that 
$x \in \fp_\pm$. This shows that $\fp_{\pm}$ is self-normalizing. 

\nin (iv) As $\g^{\pm 1}$ is abelian, the exponential 
function $\exp \colon \g^{\pm 1} \to G^{\pm 1}$ is a surjective group 
homomorphism. That it actually is a diffeomorphism follows from 
$\Ad(\exp x) h = h + [x,h] = h \mp x$ for $x \in \g^{\pm 1}$. 

\nin (v) The direct product group $G^1 \times P_-$ acts smoothly on $G$ 
by $(g,p).x := gxp^{-1}$. 
Now the orbit map $F \colon G^1 \times P_- \to G, (g,p).e = gp^{-1}$ 
has in $(e,e)$ the surjective differential 
\[ \g^1 \times \fp_- \to \g, \quad (x,y) \mapsto x - y.\] 
This implies that $F(G^1 \times P_-) = G^1 G^0 G^{-1}$ is an open subset 
of $G$ and that $F$ is a local diffeomorphism. 
By (ii), the stabilizer group of $e$ is isomorphic to 
$G^1 \cap P_- = G^1 \cap G^0 = \{e\}$. 
This follows from (ii) and the fact that, 
for $x \in \g^1$, we have $e^{\ad x} h = h + [x,h] = h - x$. 
This implies that $F$ is injective, hence a diffeomorphism onto 
an open subset of~$G$. 
\end{prf}

\begin{prop} \mlabel{prop:1.26} $P_+$ coincides 
with the flag stabilizer 
\[ G_{(\fp_+, \g^1)} := \{ g \in G \colon \Ad(g)\fp_+ = \fp_+, 
\Ad(g) \g^1 = \g^1\}.\] 
\end{prop}

\begin{prf} Clearly, $P_+ \subeq G_{(\fp_+, \g^1)}$. 
If, conversely, $g \in G_{(\fp_+, \g^1)}$, then 
$h' := \Ad(g)h$ defines a $3$-grading with 
$\g^1(h') = \g^1$ and $\g^0(h') + \g^1 = \fp_+$. 
Therefore \cite[Cor.~1.7]{BN04} implies the existence of 
$x_1 \in \g^1$ with $\ad(h' - h) = \ad x_1$, i.e., 
$h' \in h + x_1 + \fz(\g)$. Then 
$e^{\ad x_1} h' = h' + [x_1, h'] = h' - x_1 \in h + \fz(\g)$ 
shows that $\exp(x_1) g$ fixes $\ad h$, hence preserves the 
$3$-grading. This shows that 
$G_{(\fp_+, \g^1)} = G^1 G_{\ad h}$, 
so that the assertion follows from 
Lemma~\ref{lem:2grps}. 
\end{prf}

\subsection{The Olshanski semigroup $S(C_\fq)$}
\mlabel{subsec:2.3}

In this subsection we turn to the subsemigroups of $G$ 
determined by the invariant cone~$C$, resp., 
its intersections $C_\pm = \pm C \cap \g^{\pm 1}$. 
As before, $G$ is a $1$-connected Lie group with Lie algebra~$\g$,
which is $3$-graded by $\ad h$. 

\begin{prop} \mlabel{prop:1.3} 
If $C$ is pointed, then the following assertions hold: 
  \begin{itemize}
  \item[\rm(i)] The cone $C_\fq := C_+ \oplus C_- \subeq \fq = \g^1 \oplus \g^{-1}$ 
is weakly hyperbolic, i.e., $\Spec(\ad x) \subeq \R$ for $x \in C_\fq$. 
  \item[\rm(ii)] $S(C_\fq) := G^0 \exp(C_\fq)$ is a closed subsemigroup of $G$ 
invariant under $s \mapsto s^\sharp = \tau(s)^{-1}$. 
  \item[\rm(iii)] The polar map  
$\Phi \colon G^0 \times C_\fq \to S(C_\fq), (g,x) \mapsto g \exp x$ 
is a homeomorphism. 
  \end{itemize}
\end{prop}

\begin{prf} (i) From \cite[Prop.~VII.3.4]{Ne99} it follows that 
$C$ is weakly elliptic in the ideal 
$\g_C := C - C \trile \g$, i.e., 
$\Spec(\ad x) \subeq i \R$ for $x \in C$. 
For $x \in \g_C$, we have 
$\ad x(\g) \subeq \g_C$, so that 
$\Spec(\ad_\g x) \subeq \Spec(\ad_{\g_C} x) \cup\{0\}$, and
therefore $C$ is also weakly elliptic in $\g$. 

Consider the isomorphism 
\begin{equation}
  \label{eq:zeta}
 \zeta := e^{\frac{\pi i}{2}\ad h} \colon \g \to \g^c := \fh + i \fq, \quad 
x_1 + x_0 + x_{-1} \mapsto   i x_1 + x_0 - i x_{-1}.   
\end{equation}
Then $\zeta(C_\fq) = i (C \cap \fq)$ is weakly hyperbolic because 
$C \cap \fq$ is weakly elliptic.  
As $\zeta$ is an isomorphism of real Lie algebras, the cone 
$C_\fq$ is also weakly hyperbolic. 

\nin (ii), (iii): We consider the quotient group $G_{\rm ad} := \Ad(G)$ with 
Lie algebra $\g_{\rm ad} := \ad \g \cong \g/\z(\g)$ and the subgroup 
\[ G_{\rm ad}^\tau := (G_{\rm ad})^\tau 
= \{ g \in G_{\rm ad} \colon g \circ \tau = \tau \circ g \}. \] 
As $Z(G) \subeq G^0$, the Lie algebra $\g_{\rm ad}$ inherits 
a natural $3$-grading induced by $\ad h \in \g_{\rm ad}^0$. 
For every $g \in G^0$, the automorphism $\Ad(g)$ commutes with 
$\ad h$, hence also with $\tau$, so that 
$\Ad(G^0) \subeq (G_{\rm ad})^\tau$.

For the weakly hyperbolic cone $\ad(C_\fq) \subeq \ad \fq$ the relation 
$\fz(\g_{\rm ad}) \subeq \g_{\rm ad}^0$ implies that 
\[ \ad(C_\fq  - C_\fq) \cap \fz(\g_{\rm ad}) =  \{0\}.\] 
Therefore Lawson's Theorem (\cite[Thm.~7.34]{HN93}) implies that the map 
\[  \Psi_1 \colon (G_{\rm ad})^\tau \times C_\fq \to (G_{\rm ad})^\tau \exp(\ad C_\fq), 
\quad (g,x) \mapsto g \exp(x) \] 
is a homeomorphism onto a closed subset of $G_{\rm ad}$. 
Restricting to the open, hence closed subgroup $\Ad(G^0)$ of $G_{\rm ad}^\tau$, 
it follows that 
\[  \Psi \colon \Ad(G^0) \times C_\fq \to \Ad(G^0) \exp(\ad C_\fq) 
= \Ad(S(C_\fq)), \quad (g,x) \mapsto g \exp(x) \] 
is a homeomorphism onto a closed subset of $G_{\rm ad}$. 
As $Z(G) \subeq G^0$, the set $S(C_\fq) = G^0 \exp(C_\fq)$ satisfies $S(C_\fq) = S(C_\fq)Z(G)$, so that 
$S(C_\fq) = \Ad^{-1}(\Ad(S(C_\fq)))$ is a closed subset of $G$. 

Clearly, the map $\Phi$ is continuous and surjective. 
Further, the map 
\[ \sigma \colon S(C_\fq) \to C_\fq, \quad g \exp x \mapsto  x \] 
is well-defined and 
continuous because $\Psi$ is a homeomorphism and $\ad\res_\fq$ is injective 
since $\ker(\ad) = \fz(\g) \subeq \g^0$. 
Now $\Phi(g_1, x_1) = \Phi(g_2, x_2)$ implies 
$x_1 = \sigma(g_1 \exp(x_1)) =  \sigma(g_2 \exp(x_2)) = x_2$ 
which further entails that $g_1 = g_2$. Therefore $\Phi$ is injective, 
and its inverse map 
\[  \Phi^{-1} \colon S(C_\fq) \to G^0 \times C_\fq, \quad 
s \mapsto (s \cdot \exp(-\sigma(x)), x) \] 
is also continuous. This shows that $\Phi$ is a homeomorphism onto~$S(C_\fq)$. 

Now we show that $S(C_\fq)$ is a subsemigroup of $G$. 
As $G$ is $1$-connected, \cite[Cor.~7.35]{HN93} implies that 
$S(C_\fq)_0 := G^\tau_0 \exp(C_\fq)$ is a subsemigroup of $G$. 
As $\Ad(G^0)C_\fq = C_\fq$, the subgroup $G^0$ normalizes the subsemigroup 
$S(C_\fq)_0$, and thus $S(C_\fq) = G^0 S(C_\fq)_0$ also is a subsemigroup of~$G$. 

For $s \in \exp(C_\fq)$, we have $s^\sharp = s$, and for 
$s \in G^0$ we have $s^\sharp = \tau(s)^{-1} \in G^0$, so that 
$S(C_\fq)$ is $\sharp$-invariant. This completes the proof of (ii) and (iii). 
\end{prf}

\subsection{The ordered symmetric space $\cO_h$ 
and the semigroup $S(h,C)$} 
\mlabel{subsec:2.2}

The pointed closed convex cone $C \subeq \g$ defines a partial order on $\g$  by 
\[ x \leq_{C} y\quad \mbox{ if } \quad y - x \in C.\]
The invariance of $C$ under the adjoint group $\Ad(G)$ implies that 
$x \leq_C y$ implies $\Ad(g)x \leq_C \Ad(g)y$ for every $g \in G$, 
so that $G$ acts on $\g$ by order isomorphisms. 

For the formulation of the following proposition, we recall 
the concept of a symmetric space in the sense of O.~Loos: 

\begin{defn}\mlabel{defn:loos-space} Let $M$ be a smooth manifold and 
$M \times M \to M, (x,y) \mapsto x \bullet y =: s_x(y)$ 
be a smooth map with the following properties: 
each $s_x$ is an involution for which $x$ is an  isolated fixed point and 
\[ s_x(y \bullet z) = s_x(y)\bullet s_x(z) \quad \mbox{ for all } \quad x,y \in M, \quad 
\mbox{ i.e.,} \quad  s_x \in \Aut(M,\bullet).\] 
Then we call $(M,\bullet)$ a {\it symmetric space} 
(in the sense of Loos; see \cite{Lo69}). 
\end{defn}

\begin{prop}
  \mlabel{prop:1.1} {\rm(The ordered symmetric space $\cO_h$)} 
We consider the adjoint orbit $\cO_h := \Ad(G)h \cong G/G^0$. 
Then the following assertions hold: 
\begin{itemize}
\item[\rm(i)] $\cO_h$ carries the structure of a Loos symmetric space 
{\rm(Definition~\ref{defn:loos-space})}, defined by 
\[  x \bullet y = e^{\pi i \ad x} y,  \quad \mbox{resp.} \quad 
(\Ad(g)h) \bullet y = \Ad(g)\tau \Ad(g)^{-1}y.\] 
\item[\rm(ii)] $\cO_h$ carries an $\Ad(G)$-invariant 
partial order defined by restriction of $\leq_{C}$. 
\item[\rm(iii)] The order intervals $[x,y]$ in $(\cO_h, \leq)$ are compact. 
\end{itemize}
\end{prop}

\begin{prf} (i) follows from the fact that, for every 
$x \in \cO_h$, the automorphism $\tau_x := e^{\pi i \ad x}$ 
is an involutive automorphism of $\g$ for which the fixed point 
$x$ is isolated in $\cO_h =\cO_x$. 

\nin (ii) is trivial. 

\nin (iii) First we observe that the centralizer $\g^0$ of the $\ad$-diagonalizable 
element $h$ contains a Cartan subalgebra of $\g$ 
because 
every $\ad$-semisimple element is contained in a Cartan subalgebra 
(\cite[Ch.~VII, \S 2, no.~3, Prop.~10]{Bou90}). 
Now \cite[Thm.~I.13]{Ne94} implies that the 
adjoint orbit $\cO_h = \Ad(G)h \subeq \g$ is closed. 
As $\cO_h$ is closed in $\g$, the compactness of the order interval
\[ \uparrow x \cap \downarrow y = \cO_h \cap (x + C) \cap (y- C), \quad \mbox{ 
where } \quad 
\uparrow x = \{ z \in \cO_h \colon x \leq z\}, 
\downarrow y = \{ z \in \cO_h \colon z \leq y\}, \] 
follows from the compactness of $(x + C) \cap (y- C)$, a consequence 
of the pointedness of~$C$. Therefore the order intervals in 
$(\cO_h, \leq)$ are compact. 
\end{prf}

\begin{rem}
(a) Note that $\Ad(G^{\pm 1})h = h + \g^{\pm 1} \subeq \cO_h$ are affine subspaces 
of $\cO_h$ intersecting in $h$. As these subspaces are invariant under the 
stabilizer group $G^0$, they define two $G$-invariant families 
of affine subspaces of $\cO_h$. For $x = \Ad(g)h\in \cO_h$, the corresponding two 
affine subspaces through $x$ are given by 
\[ x + \g^{\pm 1}(x) = \Ad(g)(h + \g^{\pm 1}(h)).\] 

\nin (b) The affine subspace $h + \g^{\pm 1} \subeq \cO_h$ are symmetric 
subspace with respect to the canonical symmetric space structure 
defined by $h' \bullet (h' + x) = h' - x$ for $h' \in h + \g^{\pm 1}$. 
\end{rem}


Now we turn to the analysis of the semigroup 
$S(h,C) := \{ g \in G \colon \Ad(g)h \leq_C h \}$. 
We first take a closer look at an important example. 

\begin{ex} \mlabel{ex:sl2}
A typical example arises for 
\[\g = \fsl_2(\R) = \{ x \in \gl_2(\R) \colon \tr x = 0\} 
\quad \mbox{ and } \quad G = \tilde\SL_2(\R).\] 
The determinant defines an invariant Lorentzian 
form on the $3$-dimensional Lie algebra $\fsl_2(\R)$ 
and $\Ad(G)$ can be identified with the connected group $\SO_{1,2}(\R)_0$ 
acting on $3$-dimensional Minkowski space. 
Accordingly, 
\begin{align*}
C &:=
\Big \{ \pmat{ a & b \\ c & -a} \in \fsl_2(\R) \colon 
b \geq 0, c \leq 0, a^2 \leq -bc\Big\} \\
&= \Big \{ \pmat{ a & b \\ c & -a} \in \fsl_2(\R) \colon 
\det\pmat{ a & b \\ c & -a} = -a^2 - bc \geq 0, b-c \geq 0\Big\}
\end{align*}
is a pointed generating invariant closed convex cone; 
the only other one is $-C$. 
The element $h := {\displaystyle\frac{1}{2}}  \pmat{1 & 0 \\ 0 & -1}$ 
defines a $3$-grading on $\g$ with 
\[ \g^0 = \R h, \qquad \g^1 = \R E_{12}, \quad \g^{-1} = \R E_{21},\] 
and this leads to the two one-dimensional cones 
$C_+ = \R_+ E_{12}$ and $C_- = \R_+ E_{21}$. 
For $g = \pmat{a & b \\ c & d} \in \SL_2(\R)$ 
we now have 
\[ \Ad(g) h  
= \pmat{\frac{a d + bc}{2} &  - ab \\ cd & - \frac{ad + bc}{2}}
= \pmat{\frac{1}{2} + bc  &  - ab \\ cd & - \frac{1}{2} - bc}\]
and therefore 
$h - \Ad(g) h 
= \pmat{-bc &   ab \\ 
-cd &  bc} \in C$ 
if and only if 
$ab \geq 0, cd \geq 0, (bc)^2 \leq  abcd = bc(1+bc),$ 
where the latter condition is equivalent to $bc \geq 0$. We conclude that 
\[ S(h,C) = \{ g \in \SL_2(\R) \colon  ad \geq 0, cd \geq 0, bc \geq 0\}.\] 
This semigroup contains 
\[ \SL_2(\R)_+ := \{  g \in \SL_2(\R) \colon (\forall j,k)\, g_{jk} \geq 0\} 
= \exp(C_+) \exp(\R h) \exp(C_-) \] 
as its identify component, and $S(h,C) = \SL_2(\R)_+ \dot\cup - \SL_2(\R)_+$ 
(see Theorem~\ref{thm:1.7} and Corollary~\ref{cor:1.8} below). 
\end{ex}

The following lemma records some trivial 
relations between the invariant 
cone $C$ and the $3$-grading. 
\begin{lem} \mlabel{lem:C} The following assertions hold for the subalgebras 
$\fp_\pm = \g^{\pm 1} \rtimes \g^0$:  
  \begin{itemize}
  \item[\rm(i)] $C \subeq C_+ \oplus \g^0 \oplus - C_-$. 
  \item[\rm(ii)] $C_+ \oplus \g^0 \oplus C_- = \{ x \in \g \colon [h,x] \in C\} 
= (\ad h)^{-1}(C)$. 
  \item[\rm(iii)] $\downarrow h = \cO_h \cap (h - C) 
\subeq  -C_+ \oplus \g^0 \oplus  C_-$. 
  \item[\rm(iv)] $(C + \fp_\pm)/\fp_\pm = \mp C_\mp$ in $\g^{\pm 1} \cong \g/\fp_\mp$.  
  \end{itemize}
\end{lem}

\begin{prf} (i) As $C$ is invariant under $e^{\R \ad h}$, for 
$x = x_1 + x_0 + x_{-1} \in C$, we have 
\[ x_{\pm 1} = \lim_{t \to \infty} e^{- t} e^{\pm t \ad h} x \in C \cap \g^{\pm 1} 
= \pm C_\pm.\]
  
\nin (ii), (iii) follow from (i). 

\nin (iv) follows from $\mp C_\mp \subeq C$ and (i). 
\end{prf}

\begin{rem} The invariant cone $C$ is in general not uniquely determined 
by $C_+$ and $C_-$.  However, given $C_\pm$, the 
closed convex invariant cone 
\[ C^\sharp := \bigcap_{g \in G} \Ad(g)(C_+ \oplus \g^0 \oplus - C_-)\] 
contains $C$ (Lemma~\ref{lem:C}(i))  
and is contained in the product set 
$C_+ \oplus \g^0 \oplus -C_-$. It is the maximal 
invariant cone with this property. In particular, we have 
$C^\sharp \cap \g^{\pm 1} = \pm C_\pm$. 
As $C_\pm$ are pointed, the subspace 
\[ C^\sharp \cap - C^\sharp = \bigcap_{g \in G} \Ad(g)\g^0 \] 
is the largest ideal of $\g$ contained in $\g^0$. 

On the other hand, the closed convex cone $C^\flat$ generated by 
$\Ad(G)(C_+ - C_-)$ is the minimal invariant cone with 
$C^\flat_\pm = C_\pm$. 
\end{rem}

We are now ready to take a closer look at the semigroup 
$S(h,C)$ defined by the order structure on $\cO_h$. This will later 
be complemented by the result that $S(h,C) = S(C_\fq)$ 
(Theorem~\ref{thm:1.7}). 
\begin{prop} \mlabel{prop:1.5}
The set 
\[ S(h,C) = \{ g \in G \colon \Ad(g)h \leq_C h \} \] 
is a closed subsemigroup of $G$ with the following properties: 
\begin{itemize}
\item[\rm(i)] $G^0$ is the unit group $S(h,C) \cap S(h,C)^{-1}$ of~$S(h,C)$. 
\item[\rm(ii)] 
$\L(S(h,C)) := \{ x \in \g \colon \exp(\R_+ x) \subeq S(h,C) \}$ 
equals $C_+ \oplus \g^0\oplus C_-$. 
\item[\rm(iii)] $S(h,C)$ is $\sharp$-invariant. 
\item[\rm(iv)] $S(h,C) \cap  (G^{-1} G^0 G^1) = \exp(C_-) G^0 \exp(C_+) 
\subeq S(C_\fq) \subeq S(h,C)$. 
\end{itemize}
\end{prop}

\begin{prf} (i) That $S(h,C)$ is a subsemigroup 
follows immediately from the $\Ad(G)$-invariance of the order $\leq_C$, 
and its closedness follows from the closedness of $h - C$. 
Clearly $G^0$ is a subgroup of the monoid $S(h,C)$. Conversely, 
$g \in S(h,C) \cap S(h,C)^{-1}$ implies that 
\begin{equation}
  \label{eq:ginv}
\Ad(g^{-1})h - h = \Ad(g^{-1})(h - \Ad(g)h) \in 
\Ad(g^{-1})C = C 
\end{equation}
and $\Ad(g^{-1})h - h \in -C$, so that $C \cap - C = \{0\}$ yields 
$\Ad(g) h = h$, i.e., $g \in G^0$. 

\nin (ii) 
For $x \in \g^{\pm 1}(h)$ we have 
$\Ad(\exp x)h - h = e^{\ad x}h - h = [x,h] = \mp x,$ 
so that $\exp(C_+), \exp(C_-) \subeq S(h,C)$, and we thus obtain 
\begin{equation}
  \label{eq:firstincl}
\exp(C_+) G^0 \exp(C_-) \subeq S(h,C).
\end{equation}
This shows that $C_+ \oplus \g^0 \oplus C_-\subeq \L(S(h,C))$. 
Conversely, for any $x \in \L(S(h,C))$ and $t \in \R$, 
we have $e^{t \ad x}h - h \in - C$ for $t \geq 0$, 
so that 
\[ [h,x] = \lim_{t \to 0+} \frac{1}{t} (h- e^{t \ad x}h)  \in C.\] 
Now Lemma~\ref{lem:C}(ii) implies that $x \in C_+ \oplus \g^0 \oplus C_-$. 

\nin (iii) For $g \in S(h,C)$ we have 
$\Ad(g^{-1}) h - h \in C$ by \eqref{eq:ginv}, and thus 
\[ \Ad(g^\sharp) h 
= \Ad(\tau_G(g)^{-1})h 
= \tau \Ad(g)^{-1} h \in \tau(h + C) = h + \tau(C) = h - C.\]

\nin (iv) As $S(h,C)$ contains $G^0$ and 
$C_\pm \subeq \L(S(h,C))$ by (ii), we have 
\begin{equation}
  \label{eq:secincl}
\exp(C_-) G^0 \exp(C_+) \subeq S(h,C) \cap (G^{-1} G^0 G^1).
\end{equation}
As $S(C_\fq)$ is a subsemigroup, we also have 
$\exp(C_-) G^0 \exp(C_+) \subeq S(C_\fq)$. 
Further, $G^0 \subeq S(h,C)$,
 and $C_\fq = C_+ + C_- \subeq \L(S(h,C))$ yield 
$S(C_\fq) = G^0 \exp(C_\fq) \subeq S(h,C)$. 
It remains to verify that 
$S(h,C) \cap (G^{-1} G^0 G^1) \subeq \exp(C_-) G^0 \exp(C_+)$. 
So let $g = \exp(x_{-1}) g_0 \exp(x_1)$ with 
$g_0 \in G^0$ and $x_{\pm 1} \in \g^{\pm 1}$ 
and assume that $g \in S(h,C)$. Then 
\[ \Ad(g) h  \in h - C \subeq - C_+ + \fp_- \]
by Lemma~\ref{lem:C}(i).  From $[x_{-1}, \g] \subeq \fp_-$ we derive 
that $\Ad(\exp(x_{-1}))$ acts trivially on the quotient space $\g/\fp_-$. 
This leads to 
\[ - C_+ + \fp_- \ni \Ad(g_0) e^{\ad x_1} h = \Ad(g_0)(h - x_1) 
= h - \Ad(g_0) x_1,\] 
which implies that $x_1 \in \Ad(g_0)^{-1} C_+ = C_+$. 
We likewise obtain from 
\[ \Ad(g)^{-1} h  \in h + C \subeq - C_- + \fp_+ \]
that 
\[  - C_- + \fp_+ \ni \Ad(g_0)^{-1} e^{-\ad x_{-1}} h = \Ad(g_0)^{-1}(h - x_{-1}) 
= h - \Ad(g_0)^{-1} x_{-1},\] 
which implies that $x_{-1} \in \Ad(g_0) C_- = C_-$. 
Putting everything together, we see that 
\[ S(h,C) \cap  (G^{-1} G^0 G^1) \subeq \exp(C_-) G^0 \exp(C_+).\qedhere\] 
\end{prf}


To obtain finer information on $S(h,C)$, we shall 
use the Levi decomposition 
of $G$ (Lemma~\ref{lem:tridec-solrad})  
to reduce matters to the case of simple Lie algebras 
which we consider next. If $\g$ is simple, then 
$\g^1$ carries the structure of a simple euclidean Jordan algebra, 
which provides an important unifying perspective. 
For more on euclidean Jordan algebras we refer to \cite{FK94}. 
  
\begin{lem}
  \mlabel{lem:hermitian} 
If $\g$ is simple and $3$-graded by $\ad h$, then 
$\g$ is hermitian, $\g^1$ carries the structure of a euclidean Jordan 
algebra $E$, and $\g$ is isomorphic to the Lie algebra 
$\conf(E)$ of conformal vector fields on~$E$. 
For any connected Lie group $G$ with Lie algebra $\g$ 
and a maximal proper invariant cone $C \subeq \g$, we have 
\begin{equation}
  \label{eq:hermx}
S(h,C) = S(C_\fq) \subeq G^{-1} G^0 G^1.
\end{equation}
\end{lem}

\begin{prf} The first assertion $\g \cong \conf(E)$ 
follows from (\cite[Thms.~1.3.11, 3.2.8]{HO97}). 
Here $E \cong \g^1$ corresponds to the constant vector fields 
on $E$, $\g^0$ consists of linear vector fields, and 
$\g^{-1}$ of homogeneous quadratic ones. The flows of these 
vector fields generate the group $\Conf(E)_0$ 
(the identity component of the conformal group of $E$) 
which acts on $E$ by birational maps. 
Choosing a Jordan identity in $C \cap \g^1$, it follows from 
\cite[Rem.~V.4]{HNO94}  that 
$C_+ = E \cap C = \g^1 \cap C$ coincides with the positive cone 
$E_+$ of squares in the Jordan algebra~$E$. 

We first consider the $1$-connected Lie group $G$ 
with Lie algebra $\g$. Then $\Ad(G) \cong \Conf(E)_0$, 
so that we may consider the adjoint representation 
as a homomorphism $\Ad \colon G \to \Conf(E)_0$. 
In \cite[Thm.~A.1]{Ne18} we have shown that the subsemigroup 
\[ \Comp(E_+) := \{ g \in \Conf(E)_0 \colon g E_+ \subeq E_+ \} \] 
is maximal. 
We now show that $\Comp(E_+) = \Ad(S(C_\fq))$. 
In view of the polar decomposition 
\begin{equation}
  \label{eq:hermy}
\Comp(E_+) = (\Aut(E_+) \cap \Conf(E)_0)\exp(C_\fq),
\end{equation}
it suffices to show that $\Ad(G^0) = \Aut(E_+) \cap \Conf(E)_0$. 

Clearly, $\Ad(G^0)$ acts on $E \cong \g^1$ by linear maps preserving 
the positive cone $E_+ = C_+$ in the Jordan algebra~$E$. 
Suppose, conversely, that a linear automorphism $\phi$ of the convex cone 
$E_+$ is contained in the connected conformal group $\Conf(E)_0 \cong \Ad(G)$. 
Then $\phi$ defines a linear automorphism of $E$, hence 
fixes the linear vector field corresponding to $\ad h\res_{E} = \id_E$ 
(the Euler vector field on $E$). 
This means that $g \in \Ad(G)_h = \Ad(G^0)$ (Lemma~\ref{lem:2grps}). 
As explained above,  we conclude that $\Ad(S(C_\fq)) = \Comp(E_+)$. 
This implies that $\Ad(S(C_\fq))$ is a maximal subsemigroup of $\Ad(G)$. 
As its inverse image $S(C_\fq)$ in $G$ 
contains $Z(G) = \ker(\Ad)$, 
it is maximal as well. 

For the corresponding grading element $h$ and 
the maximal invariant cone $C \subeq \g$ containing $C_+$, 
this implies that the semigroup 
$S(h, C)$, which contains $S(C_\fq)$ by Proposition~\ref{prop:1.5}(iv),   
actually coincides with $S(C_\fq)$. 
Further, $S(C_\fq) \subeq G^{-1} G^0 G^1$ follows from 
Koufany's Theorem (\cite{Ko95} and also \cite[Thm.~3.8]{Ne18}). 

Any connected Lie group with Lie algebra $\g$ is of the form 
$G_\Gamma := G/\Gamma$, where $\Gamma \subeq Z(G)$ is a discrete subgroup. 
Since all three sets in \eqref{eq:hermx} are $\Gamma$-saturated, we obtain 
\[S(h,C)/\Gamma =  \{ g \in G_\Gamma \colon \Ad(g)h - h \in - C \} 
\cong S(C_\fq)/\Gamma = G_\Gamma^0 \exp(C_\fq) \subeq G_\Gamma^{-1} G_\Gamma^0 
G_\Gamma^1,\] 
and this proves \eqref{eq:hermx} for the general case.
\end{prf}

\begin{rem} In general the subgroup 
\[ \Aut(E_+) \cap \Ad(G) 
= \Aut(E_+) \cap \Conf(E)_0 \] 
is not connected. 

For $E = \Sym_n(\R)$ and $\conf(E) \cong \sp_{2n}(\R)$, 
we have for $G = \Sp_{2n}(\R)$ (not simply connected), 
$G^0 = \GL_{n}(\R)$, acting on $E$ by $g.A = g A g^\top$. 
Therefore $g = -\1$ acts trivially, so that 
\[ \Conf(E)_0  = \Ad(G) \cong \Sp_{2n}(\R)/\{\pm \1\}.\] 
If $n$ is even, then $\det(-\1) = 1$, so that 
$\GL_{n}(\R)/\{\pm \1\} \subeq \Aut(E_+)$ has two connected components. 
\end{rem}

The following theorem is the main result of this section. 
It shows that the two semigroups 
$S(C_\fq)$ and $S(h,C)$ actually coincide and 
decompose according to the $3$-grading. 

\begin{thm} \mlabel{thm:1.7} 
{\rm(Decomposition Theorem)} 
\[ S(h,C) 
= \exp(C_-) G^0 \exp(C_+) = \exp(C_+) G^0 \exp(C_-) = S(C_\fq).\]
\end{thm}

\begin{prf}  {\bf Claim 1:} $S(h,C) \subeq G^{-1} G^0 G^1$. 
In view of Proposition~\ref{prop:1.5}, 
the first equality follows from Claim~1. 
This further implies that $\exp(C_-) G^0 \exp(C_+)$ is a closed 
subsemigroup, hence coincides with the subsemigroup generated by 
$G^0$ and $\exp(C_\pm)$, and this in turn coincides with 
$S(C_\fq) = G^0 \exp(C_\fq)$. We also obtain from Proposition~\ref{prop:1.3}  that 
\[ S(h,C) = S(h,C)^\sharp = ( \exp(C_-) G^0 \exp(C_+))^\sharp 
= \exp(C_+) G^0 \exp(C_-).\] 
So it remains to verify Claim~1. 

In view Lemma~\ref{lem:tridec-solrad}, $\g$ contains 
an $\ad h$-invariant Levi complement~$\fl$. 
Let $\fl = \fl_1 \oplus \cdots \oplus \fl_m$ denote the decomposition 
into simple ideals. Then each ideal is invariant under $\ad h$ 
because all derivations of $\fl$ are inner. 
If $\fl_j$ is compact, then $\fl_j \subeq \fl^0$ because all derivations of 
$\fl_j$ are elliptic. If the grading of $\fl_j$ is non-trivial, 
then $\fl_j$ is contained in the ideal 
$\g_C := C - C$ of $\g$ which contains $\g^{\pm 1} = C_\pm - C_\pm$. 
By Lemma~\ref{lem:hermitian}, $\fl_j$ is hermitian. 
It follows in particular that $\fl^1 \oplus \fl^{-1}$ is contained 
in the sum $\fl_h$ of all simple hermitian ideals of $\fl$.

\nin {\bf Claim~2:} $\oline{p_{\fl_h}(C)} \subeq \fl_h$ 
is a pointed generating invariant cone. 
In $\g_C$ the cone $C$ is pointed and generating, 
so that $\g_C$ contains a compactly embedded Cartan subalgebra 
$\ft$, compatible with the Levi decomposition $\g_C = \fr \rtimes \fl$
(\cite[Prop.~VII.1.9]{Ne99}). By \cite[Thm.~VII.3.8]{Ne99} 
there exists an adapted positive system $\Delta^+ \subeq \Delta(\g_C,\ft)$ 
(cf.\ Appendix~\ref{app:c}) such that 
\[ C \cap \ft \subeq C_{\rm max}(\Delta_p^+) 
= (i\Delta_p^+)^\star.\] 
Here $\Delta_p \subeq \Delta$ is the subset of non-compact roots 
which contains the subset $\Delta_{p,s}$ of non-compact simple 
roots corresponding to $\fl_h$. All these roots vanish on $\ft \cap  \fr$, so that 
\[ (i\Delta_p^+)^\star \subeq 
(i\Delta_{p,s}^+)^\star = (i\Delta_{p,s}^+)^\star + \ft \cap \fr. \] 
This implies that 
\[ p_{\fl_h}(C \cap \ft) \subeq (i\Delta_{p,s}^+)^\star,\] 
and since $\Ad(G)(C\cap \ft)$ is dense in $C$, it follows that 
\begin{equation}
  \label{eq:dag1}
 p_{\fl_h}(C) \subeq \oline{\Ad(L_h) C_{\rm max}(\Delta_{p,s}^+)} 
\subeq W_{\rm max,s} := \oline{\Ad(L)(C_{\rm max}(\Delta_{p,s}^+))}.   
\end{equation}
Here $W_{\rm max,s}$ is a pointed closed convex invariant cone in $\fl_h$. 
As $C$ generates $\g_C$, the cone $p_{\fl_h}(C)$ generates $\fl_h$. 
This proves Claim~2. 

If the grading on a hermitian ideal $\fl_j$ is non-trivial, 
the projection of $C_+$ into $\fl_j^1$ is contained 
in a pointed invariant cone, and this in turn implies that 
$\fl_j^1$ can be identified with a euclidean Jordan algebra $E_j$ 
for which $\fl_j\cong \conf(E_j)$ is the Lie algebra of the conformal group. 
This follows from Lemma~\ref{lem:hermitian}, which further entails that, 
for the invariant cone $C_j := W_{\rm max,s} \cap \fl_j$ in $\fl_j$, 
we have 
\[  S(h_j, C_j)  = S(C_{j,\fq}) \subeq L_j^{-1} L_j^0 L_j^1.\] 
We conclude that, for the grading element $h_\fl = \sum_{j=1}^m h_j$ 
of the $3$-graded semisimple Lie algebra $\fl$, we have  
\begin{equation}
  \label{eq:dag2}
S\Big(h_\fl, \sum_j C_j\Big) = \prod_{j = 1}^m S(h_j, C_j)  \subeq L^{-1} L^0 L^{1}, 
\end{equation}
where we use that $L_j = L_j^0$ if $h_j = 0$. 

If $C_1 \subeq C_2$ are invariant cones in $\fl$, then we clearly have 
$S(h_\fl, C_1) \subeq S(h_\fl, C_2),$ so that \eqref{eq:dag1} and \eqref{eq:dag2} 
show that 
\begin{equation}
  \label{eq:dag3}
 p_L(S(h,C)) \subeq S(h_\fl, W_{\rm max,s}) 
\subeq L^{-1} L^0 L^{1}.
\end{equation}
With Lemma~\ref{lem:tridec-solrad}(iii), we now obtain 
$S(h,C) \subeq p_L^{-1}(L^{-1}L^0 L^1) = G^{-1} G^0 G^1.$ 
\end{prf}

The subgroup $G^\tau$ is connected because $G$ is $1$-connected 
(\cite[Thm.~B.2]{Ne19}), so 
it coincides with $(G^0)_0$, and the preceding theorem implies that: 

\begin{cor}\mlabel{cor:1.8} The identity component of $S(C_\fq)$ is 
\[ S(C_\fq)_0 = G^\tau \exp(C_\fq) = \exp(C_+) G^\tau\exp(C_-).\]
\end{cor}

\begin{rem} The Decomposition Theorem shows in particular that the semigroup 
$S(h,C)$ only depends on the cone $C_\fq = C_+ \oplus C_-$, i.e., 
that $S(h,C) = S(h,C')$ if $C_\pm = C_\pm'$. 
As we have seen in the proof of Lemma~\ref{lem:hermitian}, this only leads to 
two different semigroups $S(h,C)$ and $S(h,-C) = S(h,C)^{-1}$. 
\end{rem}

\begin{ex} (a) Suppose that $\g$ is solvable and that 
$C_\pm$ span $\g^{\pm 1}$. 
Then $[\g,\g]$ is a nilpotent ideal containing 
$\g^{\pm 1}$. For any pointed invariant cone $C$, the cone 
$C \cap [\g,\g]$ is a pointed invariant 
cone in the nilpotent Lie algebra $[\g,\g]$, so that its span 
is abelian by \cite[Ex.~VII.3.21]{Ne99}. Then 
$[\g^1, \g^{-1}] = \{0\}$ and 
\[ \g \cong (\g^1 \oplus \g^{-1}) \rtimes \g^0.\] 

Conversely, any involution $D \colon E \to E$ of a finite dimensional 
vector space $E$ defines a solvable Lie algebra 
\[ \g := E \rtimes_D \R \quad \mbox{ with the bracket } \quad 
[(v,t), (v',t')] := (t Dv' - t' Dv, 0).\] 
For $h := (0,1)$ we then obtain the $\ad h$-eigenspaces 
\[ \g^0 = \R h \quad \mbox{ and } \quad \g^{\pm 1} = E^\pm(D).\] 

\nin (b) Let $(V,\omega)$ be a symplectic vector space and  
$\heis(V,\omega) = \R \oplus V$ be the corresponding Heisenberg algebra 
with the bracket $[(z,v),(z',v')] = (\omega(v,v'),0)$. 
Then any pointed invariant cone $C$ is contained in the center, 
hence (up to sign) of the form $C = \R_+ (1,0)$. 
If $\tau$ is an involutive automorphism of $\heis(V,\omega)$ 
with $\tau(C) = - C$, then (up to equivalence) it has the form 
\[ \tau(z,v) = (-z,\tau_V(v)),\] 
where $\tau_V \colon V \to V$ is antisymplectic, i.e., $\tau_V^*\omega = - \omega$. 

Extending $\heis(V,\omega)$ by a diagonalizable derivation 
$D$ to $\g := \heis(V,\omega) \rtimes_D \R$, we may also 
consider the corresponding element $h := (0,0,1)$ for which 
$\ad h$ coincides with $D$ on the Heisenberg algebra 
and extend $\tau$ by $\tau(h) = h$ to $\g$ 
(which works if $D$ commutes with $\tau_V$). 
Suppose that $\ad h$ defines a $3$-grading with $\tau = e^{\pi i \ad h}$ 
and, w.l.o.g., that $Dz = z$ for the central element $z = (1,0)\in 
\heis(V,\omega)$. 
Then $V_\fh := \Fix(\tau_V)$ and $V_\fq = \Fix(-\tau_V)$ are 
Lagrangian subspaces with $V = V_\fh \oplus V_\fq$ and 
$\g^0 = V_\fh \oplus \R h$. 
From $[V_\fq, V_\fh] \subeq \R z \subeq \g^1$, it follows that 
$V_\fq \subeq \g^1$. We therefore have
\[ \g^{-1} = \{0\} \quad \mbox{ and } \quad 
\g^1 = V_\fq = \Fix(-\tau_V).\] 
\end{ex}

For more complicated examples we refer to 
Subsection~\ref{subsec:exs}. 

\subsection{The semigroup $G_{\pi i}(C)$} 
\mlabel{subsec:2.6}

Let $\eta_G \colon G \to G_\C$ denote the universal complexification of $G$, 
i.e., $G_\C$ is the $1$-connected Lie group with Lie algebra 
$\g_\C$, and $\L(\eta_G) \colon \g \to \g_\C$ is the canonical 
inclusion. 
For the pointed generating invariant cone $C$, 
Lawson's Theorem (\cite[Thm.~IX.1.10]{Ne99}) implies the existence 
of a semigroup $S(iC)$ which is a covering of the subsemigroup 
$\eta_G(G) \exp(i C)$ of the universal complexification $G_\C$ 
(the simply connected group with Lie algebra $\g_\C$). 
Then the exponential function $\exp \colon \g + i C \to G_\C$ 
lifts to an exponential function $\Exp \colon \g + i C \to S(iC)$ 
and the polar map 
\[  G \times C \to S(iC) = G \Exp(iC), \quad 
(g,x) \mapsto g \Exp(ix) \] 
is a homeomorphism, and, if $C$ has non-empty interior,\begin{footnote}
{Note that we do not assume that $C$ has interior points. A typical 
example where this is not the case arises from the Poincar\'e group; 
see Example~\ref{ex:poincare}.}  
\end{footnote}
a diffeomorphism of $G \times C^0$ onto the complex manifold $S(iC)^0$. 

For $z \in \C$ with $\Im z > 0$, we write $G_z = G_z(C) 
\subeq G$ for the closed subsemigroup of all elements $g \in G$ 
for which the 
orbit map $\beta^g(t) := \exp(th) g \exp(-th)$ extends 
analytically to a continuous map 
\[ \beta^g \colon  \oline{\cS_{\Im z}} = \{ w \in \C \colon 0 \leq \Im w \leq \Im z \} \to S(iC) \]
(see \cite[Lemma~3.9]{Ne19} for details). Here ``analytic'' means that, on 
the open strip the composed map $\beta^g \colon \cS_{\Im z} \to S(iC) \to G_\C$ 
is holomorphic. 

\begin{lem}
  \mlabel{lem:abcase} 
$(G^1)_{\pi i} = \exp(C_+)$ and $(G^{-1})_{\pi i} = \exp(C_-)$. 
\end{lem}

\begin{prf} For the abelian subgroup $G^1 \cong \g^1$ and 
$G^1_\C \cong \g_\C^1$, we have 
\[ G^1_{\pi i} = \{ \exp(x) \colon x \in \g^1, (\forall y \in [0,\pi])\ 
\exp(e^{iy} x) \in G^1 \exp(i C_+) = \exp(\g_1 + i C_+) \}.\]
As the exponential function of $(G^1)_\C$ is bijective,\begin{footnote}
{This follows from the same argument as 
Lemma~\ref{lem:conefield-a}(iv).}\end{footnote}
$\exp(x) \in G^1_{\pi i}$ is equivalent to 
\[ e^{iy} x = \cos(y) x + i \sin(y) x \in \g^1 + i C_+ \quad \mbox{ for } \quad 
y \in [0,\pi],\] 
and this is equivalent to $x \in C_+$. 
This proves the first assertion, and the second follows similarly.   
\end{prf}

\begin{thm} \mlabel{thm:1.29}
$G_{\pi i} = \exp(C_-) G^0 \exp(C_+)$. 
\end{thm}

\begin{prf} 
{\bf Step 1} $S(h,C) \subeq G_{\pi i}(C)$: \\
By Theorem~\ref{thm:1.7}, we have 
$S(h,C)  = \exp(C_+) G^0 \exp(C_-).$ 
Since $G_{\pi i}(C)$ is a subsemigroup of $G$ which obviously contains 
$G^0$ (the elements with constant orbit maps), it suffices to see that 
$\exp(C_{\pm}) \subeq G_{\pi i}$. This follows from 
Lemma~\ref{lem:abcase}. 

\nin {\bf Step 2.} 
 $G_{\pi i}(C) \subeq G^{-1} G^0 G^1$: \\
Let $G \cong R \rtimes L$ be a Levi decomposition with 
$[h,\fl] \subeq \fl$ (Lemma~\ref{lem:tridec-solrad}) 
and write $p_L \colon G \to L$ for the corresponding 
morphism of $3$-graded Lie groups. In view of Lemma~\ref{lem:tridec-solrad}(iii), 
it suffices to show that $p_L(G_{\pi i}) \subeq L^{-1} L^0 L^1$. 

We have already seen in the proof of Theorem~\ref{thm:1.7} 
that $C_\fl := \oline{q_\fl(C)} \subeq \fl$ 
is a pointed invariant cone whose span $\fl_C$ is a direct sum of hermitian ideals. 
All other simple ideals of $\fl$ are contained in~$\fl^0$. 
For $L_{\pi i} = L_{\pi i}(C_\fl)$, 
it follows that $p_L(G_{\pi i}) \subeq L_{\pi i}$. Enlarging the cone 
$C_\fl$ to a maximal pointed invariant cone $C_{\rm max}$ in $\fl_\C$, 
we have $L_{\pi i}(C_\fl) \subeq L_{\pi i}(C_{\rm max})$. 
As $C_{\rm max}$ is adapted to the decomposition 
$\fl_C = \fl_1 \oplus \cdots \oplus \fl_m$ into simple ideals in the sense that 
\[ C_{\rm max} = \oplus_{j = 1}^m (C_{\rm max} \cap \fl_j), 
\quad \mbox{ it follows that } \quad 
 L_{\pi i}(C_{\rm max}) = \prod_{j = 1}^m L_{j,\pi i}(C_{\rm max} \cap \fl_j). \] 

Therefore it suffices to show that $L_{\pi i} \subeq L^{-1} L^0 L^1$ if $L$ 
is simple hermitian and $C = C_{\rm max}$. Then 
$L_{\pi i}$ is a closed subsemigroup 
of $L$ containing the maximal subsemigroup 
\[ \exp(C_+) L^0 \exp(C_-) = S(C_\fq) \] (Lemma~\ref{lem:hermitian}). 
As $(L^1)_{\pi i} = L^1 \cap L_{\pi i} = \exp(C_+)$ follows 
from Lemma~\ref{lem:abcase}, $L_{\pi i}\not= L$, so that 
the maximality of $S(C_\fq)$ implies that 
$L_{\pi i} =  S(C_\fq) \subeq L^{-1} L^0 L^1.$

\nin {\bf Step 3.}  $G_{\pi i}(C) \subeq \exp(C_+) G^0 \exp(C_-) = S(h,C)$: \\ 
In view of Step 2, it remains to show that 
$g = g_1 g_0 g_{-1} \in G_{\pi i}$ with $g_j \in G^j$ implies that 
$g_{\pm 1} \in \exp(C_\pm)$. 
To this end, we consider the projections 
\[ \pi_\pm \colon G \to M_\pm := G/P_\mp.\] 
Clearly, $\pi_+$ maps $G_{\pi i}$ into the subset of all elements 
$m \in M_+$ for which the orbit map 
$\gamma^m(t) := \exp(th).m$ extends analytically to a map 
from $\oline{\cS_\pi}$ with values in 
\[ \pi_-(\exp(\g^1 + i C_+)) \subeq M_{+,\C} := G_\C/P_{-,\C}.\] 
Writing $g_1 = \exp(x)$,  we have 
$\pi_+(g) = g_1 P_- = \exp(x) P_-$ and 
\[ \gamma^{\pi_+(g)}(z) 
= \exp((a + i b)h) \exp(x) P_{-,\C} 
= \exp(e^a e^{ib} x) P_{-,\C} 
= \exp\big(e^a( \cos(b) x + i \sin(b)x)\big) P_{-,\C}.\]  
As $\eta_{+,\C} := \g^{1}_\C \to M_{+,\C}, z \mapsto \exp(z) P_-$ 
also is an open embeddings (Lemma~\ref{lem:conefield-a}(v)), we see that 
$g \in G_{\pi i}$ implies $x \in C_+$. 
A similar argument shows that $g_{-1} \in \exp(C_-)$. 

Combining Steps 1-3, the assertion follows.   
\end{prf}

\begin{cor} \mlabel{cor:suinv} 
We have 
\[ G_{\pi i}(C,\tau_G) 
:= \{ g \in G_{\pi i}(C) \colon \beta^g(\pi i) = \tau_G(g)\} 
=  \exp(C_-) (G^0)^\tau \exp(C_+).\] 
\end{cor}

\begin{prf} By Theorem~\ref{thm:1.29}, we only have to determine which 
elements in $G_{\pi i} = \exp(C_-) G^0 \exp(C_+)$ satisfy 
$\beta^g(\pi i) = \tau_G(g)$. Writing $g = g_{-1} g_0 g_1$ with 
$g_{\pm 1} \in \exp(C_\pm)$, we have $\beta^{g_{\pm 1}}(\pi i) = g_{\pm 1}^{-1} 
= \tau_G(g_{\pm 1})$. We thus obtain 
\[ \beta^g(\pi i) = g_{-1}^{-1} g_0 g_1^{-1} 
\quad \mbox{ and } \quad 
\tau_G(g) = g_{-1}^{-1} \tau_G(g_0) g_1^{-1}.\] 
Equality of these elements is equivalent to $\tau_G(g_0) = g_0$.
\end{prf}


\begin{cor} \mlabel{cor:suinv2} 
If $(U,\cH)$ is an antiunitary representation of 
$G \rtimes \{\1,\tau_G\}$ with discrete kernel, 
$J_\sV = U(\tau_G)$, $\Delta_\sV = e^{2\pi i \partial U(h)}$,  
and $C = C_U$, then 
\[ \{ g \in G_{\pi i} \colon U(\beta^g(\pi i)) = U(\tau_G(g)) \} 
= \exp(C_-) G_\sV\exp(C_+) \subeq S_\sV,\] 
where $G_\sV = \{ g \in G \colon U(g) \sV = \sV\}$. 
\end{cor}

\begin{prf} As in the proof of Corollary~\ref{cor:suinv}, 
we see that $g = g_{-1} g_0 g_1$ is contained in the set on the left hand side 
if and only if $U(g_0) = U(\tau_G(g_0)) = J_\sV U(g_0) J_\sV$. 
Since $\ker(U)$ is discrete, the assertion now follows from 
\begin{align*}
 G_\sV 
&= \{ g \in G \colon U(g) J_\sV = J_\sV U(g), U(g) \Delta_\sV U(g)^{-1} 
= \Delta_\sV \} 
= \{ g \in G^0 \colon U(g) J_\sV = J_\sV U(g) \}.
\end{align*}
The inclusion $\exp(C_-) G_\sV\exp(C_+) 
\subeq S_\sV$ now follows from the Inclusion Theorem, which is  \cite[Thm.~3.11]{Ne19}. 
\end{prf}

\section{The semigroup $S_\sV$ in the $3$-graded case} 

In this section we eventually turn to the compression 
semigroup $S_\sV$ of a standard subspace~$\sV$. So we consider an  
antiunitary representation $(U,\cH)$ of 
a semidirect product $G \rtimes \{\id_G, \tau_G\}$ with discrete 
kernel, where $G$ is a connected Lie group, 
$\tau_G$ is an involutive automorphism of $G$ and $h\in \g^\tau$. 
We consider the standard subspace  $\sV \subeq \cH$ specified by 
the modular objects 
$J_\sV = U(\tau_G)$ and $\Delta_\sV = e^{2\pi i \partial U(h)}$. 

\subsection{The general monotonicity theorem} 

The following result is essentially contained in the work of 
Borchers and Wiesbrock, see for instance \cite[\S II.1]{Bo00}. 
For the formulation we refer to the discussion of the order 
on the space of selfadjoint operators 
in Appendix~\ref{app:a}. 

\begin{thm} {\rm(Borchers--Wiesbrock Monotonicity)}  \mlabel{thm:b.2}
If $\sV_1 \subeq \sV_2$ are standard subspaces of $\cH$, then 
\[ \Delta_{\sV_2} \leq \Delta_{\sV_1},\] 
and we also have $\log(\Delta_{\sV_2}) \leq \log(\Delta_{\sV_1})$ in the sense that 
\begin{equation}
  \label{eq:logdelt}
q_{\log(\Delta_{\sV_2})}(\xi,\xi) \leq q_{\log(\Delta_{\sV_1})}(\xi,\xi) 
\quad \mbox{ for } \quad 
\xi \in \cD[\log(\Delta_{\sV_2})] \cap \cD[\log(\Delta_{\sV_1})].
\end{equation}
\end{thm}

\begin{prf} Let $T_{\sV_1} \subeq T_{\sV_2}$ be the Tomita operators 
of $\sV_1$ and $\sV_2$, respectively. Their graphs 
\[ \Gamma_j := \Gamma(T_{\sV_j}) = \{ (\xi, T_{\sV_j} \xi) \colon \xi \in \cD(T_{\sV_j}) \} \subeq 
\cH \oplus \cH \]  
are closed subspaces of $\cH \oplus \cH$ with 
$\Gamma_1 \subeq \Gamma_2$. Hence the orthogonal projections 
$P_j$ on $\cH \oplus \cH$ with range $\Gamma_j$ satisfy 
$P_1 \leq P_2$. Identifying $B(\cH \oplus \cH)$ with the algebra 
$M_2(B(\cH))$ of $(2 \times 2)$-matrices with entries in $B(\cH)$, we 
write $P_j$ as a $(2 \times 2)$-matrix. We obtain from 
$\sV_1 \subeq \sV_2$ that $\Gamma(T_{\sV_1}) \subeq \Gamma(T_{\sV_2})$, 
i.e., $P_1 \leq P_2$, and hence that $(P_1)_{11} \leq (P_2)_{11}$. 
Therefore Lemma~\ref{lem:b.1} leads to 
\[(\1 + \Delta_{\sV_1})^{-1} 
=  (\1 + T_{\sV_1}^* T_{\sV_1})^{-1} \leq (\1 + T_{\sV_2}^* T_{\sV_2})^{-1} 
= (\1 + \Delta_{\sV_2})^{-1}.\] 
As the function $x \mapsto -\frac{1}{x}$ is operator monotone on $(0,\infty)$ 
(cf.~\cite[Cor.~10.13]{Sch12}),  we obtain 
$\Delta_{\sV_2} \leq \Delta_{\sV_1},$ and \eqref{eq:logdelt} follows from 
Theorem~\ref{thm:logmon}. 
\end{prf}

\begin{rem} Note that the relation 
$\Delta_{\sV_2}\leq \Delta_{\sV_1}$ conversely implies that 
\begin{equation}
  \label{eq:delrel}
\sV_1 + i \sV_1 = \cD(\Delta_{\sV_1}^{1/2}) = \cD[\Delta_{\sV_1}]  
\subeq \cD[\Delta_{\sV_2}] = \cD(\Delta_{\sV_2}^{1/2}) = \sV_2 + i \sV_2 
\end{equation}
(Definition~\ref{def:a.2}). In general, the inclusion \eqref{eq:delrel} 
is weaker than $\Delta_{\sV_2}\subeq \Delta_{\sV_1}$. 
If, for instance, $\Delta_{\sV_1}$ and $\Delta_{\sV_2}$ are bounded, 
then $\sV_1 + i \sV_1 = \sV_2 + i \sV_2 = \cH$, 
but $\Delta_{\sV_2}\leq \Delta_{\sV_1}$ does not always hold. 
\end{rem}

We now apply Theorem~\ref{thm:b.2} to obtain information on $S_\sV$. 

\begin{thm} {\rm(The Monotonicity Theorem)}  \mlabel{thm:inc}
Let $(U,\cH)$ be an antiunitary representation of 
$G \rtimes \{\id_G, \tau_G\}$, $h \in \g^\tau$, 
and let $\sV \in \Stand(\cH)$ be determined by 
\[ J_\sV = U(\tau_G) \quad \mbox{ and } \quad 
\Delta_\sV = e^{2\pi i \partial U(h)}.\] 
Then 
\[ S_\sV = \{ g \in G \colon U(g) \sV \subeq \sV \} 
\subeq  S(h,C_U) = \{ g \in G \colon \Ad(g)h \in h - C_U\}.\] 
\end{thm}

\begin{prf} For $g \in S_\sV$, we have $U(g)\sV \subeq \sV$. 
As $\Delta_{U(g)\sV} = e^{2\pi i \partial U(\Ad(g)h)}$ and 
\[ \cH^\infty \subeq 
\cD(i\partial U(h)) \cap \cD(i\partial U(\Ad(g)h)) 
\subeq \cD[i\partial U(h)] \cap \cD[i\partial U(\Ad(g)h)],\] 
Theorem~\ref{thm:b.2} implies that 
\begin{equation}
  \label{eq:relonsmo}
\la \xi, i \partial U(\Ad(g)h)\xi \ra \geq 
\la \xi, i \partial U(h) \xi \ra \quad \mbox{ for } \quad 
\xi \in \cH^\infty.  
\end{equation}
As the operators 
$i\partial U(x)$, $x \in \g$, 
are the closures of their restriction to the 
$U(G)$-invariant subspace~$\cH^\infty$ of smooth vectors, 
we conclude from \eqref{eq:relonsmo} that $\Ad(g)h - h \in - C_U,$ 
so that $\Ad(g)h \in h - C_U$, i.e., $g \in S(h,C_U)$. 
This proves the theorem. 
\end{prf}

\subsection{The semigroup $S_\sV$ in the $3$-graded case} 

In the context that we studied throughout this paper, where 
$\g$ is $3$-graded by $\ad h$, we have the following stronger result: 

\begin{thm} \mlabel{thm:sv}
Suppose, in addition to the setting of {\rm Theorem~\ref{thm:inc}},
 that $\g$ is $3$-graded by $\ad h$, $\tau = e^{\pi i \ad h}$ 
and $C = C_U$, where $\ker(U)$ is discrete. Then 
\begin{equation}
  \label{eq:svident}
 S_\sV = \exp(C_+) G_\sV \exp(C_-). 
\end{equation}
\end{thm}

\begin{prf} First we recall from 
Theorem~\ref{thm:1.7} that, under the stated assumptions,  
$S(h,C) = \exp(C_+) G^0 \exp(C_-)$. 
Next Theorem~\ref{thm:inc} shows that $S_\sV \subeq S(h,C)$. 
With Corollary~\ref{cor:suinv2}, we thus obtain 
\begin{equation}
  \label{eq:firstinc}
 \exp(C_+) G_\sV \exp(C_-) \subeq S_\sV \subeq S(h,C) 
= \exp(C_+) G^0 \exp(C_-).
\end{equation}

Let $g = g_+ g_0 g_-$ with $g_\pm \in \exp(C_\pm)$ and $g_0 \in G^0$ be 
an element of $S(h,C)$. 
If $g \in S_\sV$, then $U(g)\sV \subeq \sV$ implies that 
the orbit map $\alpha^{U(g)}(t) := \Delta_\sV^{-it/2\pi} U(g) \Delta_\sV^{it/2\pi}$ 
extends to $\oline{\cS_\pi}$ with 
\begin{equation}
  \label{eq:az-cond}
\alpha^{U(g)}(\pi i) = U(\tau_G(g))
\end{equation}
(see the Araki--Szid\'o Theorem; \cite{AZ05}, \cite[Thm.~2.3]{Ne19}). 
We know already that $\alpha^{U(g_\pm)}(\pi i) = U(g_\pm)^{-1}$ exists, 
and, since $\alpha^{U(g_0)}$ is constant with $\alpha^{U(g_0)}(\pi i) = U(g_0)$, 
we obtain from \eqref{eq:az-cond} 
\[ \alpha^{U(g)}(\pi i) 
= \alpha^{U(g_+)}(\pi i) \alpha^{U(g_0)}(\pi i) \alpha^{U(g_-)}(\pi i) 
 = U(g_+)^{-1} U(g_0)U(g_-)^{-1} \] 
and 
\[ U(\tau_G(g)) = U(g_+)^{-1} U(\tau_G(g_0)) U(g_-)^{-1} \] 
that $U(g_0) = U(\tau_G(g_0)),$ 
so that $g_0 \in G_\sV$. This shows that 
\[ S_\sV = S_\sV \cap S(h,C) \subeq \exp(C_+) G_\sV\exp(C_-),\] 
 and with \eqref{eq:firstinc} we obtain 
\eqref{eq:svident}.
\end{prf}

\subsection{Examples} 
\mlabel{subsec:exs}

\begin{ex} \mlabel{ex:poincare} (Poincar\'e group) 
In Quantum Field Theory on Minkowski space, 
the natural symmetry group  is 
the proper Poincar\'e group 
$P(d) \cong \R^{1,d-1} \rtimes \OO_{1,d-1}(\R)^\uparrow$ 
acting by orientation preserving 
isometries on $d$-dimensional Minkowski space $\R^{1,d-1}$. 
Its Lie algebra is  $\g := \fp(d) \cong \R^{1,d-1} \rtimes \so_{1,d-1}(\R)$ 
and the closed forward light cone 
\begin{equation}
  \label{eq:lightcone}
C := \{ (x_0, \bx) \in \R^{1,d-1} \colon x_0 \geq 0, x_0^2 \geq \bx^2\} 
\end{equation}
is a pointed invariant cone in $\fp(d)$. 
The generator $h \in \so_{1,d-1}(\R)$ of the Lorentz boost on the 
$(x_0,x_1)$-plane 
\[ h(x_0,x_1, \ldots, x_{d-1}) = (x_1, x_0, 0,\ldots, 0)\] 
defines a $3$-grading on $\g$ because 
$\ad h$ is diagonalizable with spectrum $\{-1,0,1\}$, 
and $\tau := e^{\pi i \ad h}$ defines an involution on 
$\g$, acting on the ideal $\R^{1,d-1}$ (Minkowski space) by 
\[ \tau_M(x_0, x_1, \ldots, x_{d-1}) = (-x_0, -x_1, x_2, \ldots, x_{d-1}).\] 
To connect with the results above, we have to apply them 
to the universal cover $\tilde G$ of the group 
$G := P(d)_0  \cong \R^{1,d-1} \rtimes \SO_{1,d-1}(\R)^\uparrow$. 

A unitary representation $(U,\cH)$ 
of $G$ is called a {\it positive energy representation} 
if $C \subeq C_U$. If $\ker(U)$ is discrete, then $C_U$ is pointed, 
and $C = C_U$ follows from the fact that 
this is, up to sign, the only non-zero pointed
invariant cone in the Lie algebra $\g = \fp(d)$ for $d  > 2$; 
for $d = 2$ there are four pointed invariant  cones which 
are quarter planes. 

The Lie algebra $\g$ is $3$-graded by $\ad h$, 
but $\g^0$ and the two cones $C_\pm$ generate a proper Lie subalgebra. 
Here $\g^0 = \ker(\ad h) = \fh$  is the centralizer of the Lorentz boost: 
\[ \g^0 = (\{(0,0)\} \times \R^{d-2}) 
\rtimes (\so_{1,1}(\R) \oplus \so_{d-2}(\R)) 
\cong (\R^{d-2} \rtimes \so_{d-2}(\R)) \oplus \R h, \] 
and, 
\[ C_+ = C \cap \g^1 = \R_+ (e_1 + e_0)
\quad \mbox{ and } \quad 
C_- = -C \cap \g^{-1} = \R_+ (e_1 - e_0).\] 

The subsemigroup $S(h,C) := \{ g \in G \colon h - \Ad(g)h \in C\}$ 
is easy to determine. 
The relation $\Ad(g)h-h \in \R^d$ implies that $g = (v,l)$ with 
$\Ad(l)h = h$, and then 
$\Ad(g)h = \Ad(v,\1)h = -hv \in - C$ is equivalent to $hv \in C$, which 
specifies the closure $\oline{W_R}$ of the standard right wedge 
\[ W_R := \{ x \in \R^{1,d-1} \colon x_1 > |x_0|\}. \] 
We therefore obtain 
\[ S(h,C) = \oline{W_R} \rtimes \big(\SO_{1,1}(\R)^\uparrow \times \SO_{d-2}(\R)\big) 
= \{ g \in G \colon gW_R \subeq W_R\}\]
(see \cite[Lemma~4.12]{NO17} for the last equality).  
As the subgroup 
$\SO_{1,1}(\R)^\uparrow \times \SO_{d-2}(\R) \subeq \SO_{1,d-1}(\R)$ 
commutes with $h$ and $\tau$. 
For any antiunitary positive energy representation of 
\[ G \rtimes \{\1,\tau_M\} = \R^{1,d-1} \rtimes \OO_{1,d-1}(\R)^\uparrow, \] 
the semigroup $S_\sV$ corresponding to the standard subspace 
specified by $U(\tau_M) = J_\sV$ and $\Delta_\sV = e^{2\pi i \partial U(h)}$ 
satisfies 
\[ S_\sV = S(h,C) = \oline{W_R} \rtimes \big(\SO_{1,1}(\R)^\uparrow \times \SO_{d-2}(\R)\big), 
\quad \mbox{ where } \quad 
\SO_{1,1}(\R)^\uparrow = \exp(\R h).\] 
For the covering group $\tilde G$ we obtain the same picture 
because the involution acts trivially on the 
covering $(\tilde G)^0$ of~$G^0$. 
\end{ex}

\begin{ex}
  (Conformal groups $\SO_{2,d}(\R)$) 
The Lie algebra of the conformal group 
$G := \SO_{2,d}(\R)^\uparrow$ of Minkowski space is 
$\g = \so_{2,d}(\R)$, which contains the Poincar\'e algebra 
as those elements corresponding to affine vector fields on $E := \R^{1,d-1}$. 
For $d \geq 3$ it is a simple hermitian Lie algebra. 
It contains many elements $h$ defining a $3$-grading on $\g$, 
but all these elements are conjugate. One arises from the 
element $h_0 = \id$ corresponding to the Euler vector field on~$E$. 
Then $\g^j(h_0)$, $j = -1,0,1$, are spaces of vector fields 
on $E$ which are linear (for $j = 0$), constant (for $j = 1$) and 
quadratic (for $j = -1$).\begin{footnote}
{We encountered similar 
structures in the proof of Lemma~\ref{lem:hermitian} for more 
general euclidean Jordan algebras.} \end{footnote}
Another important example is the 
element $h_1 \in \so_{1,1}(\R) \subeq \so_{2,d-1}(\R)$ 
corresponding to a Lorentz boost in the Poincar\'e algebra 
(see Example~\ref{ex:poincare}). 

We consider the minimal invariant cone $C \subeq \g$ 
which intersects $E$ in the positive light cone~$C_+(h_0)\subeq E$. 
For all these elements $h$ we obtain a complete description of 
the corresponding semigroups $S_\sV$ as 
$\exp(C_+) G_\sV \exp(C_-)$, and here these semigroups 
have interior points because $C_\pm$ generate the subspaces $\g^{\pm 1}$. 
\end{ex}

\begin{ex} Another interesting example which is 
neither semisimple nor an affine group is given by 
the Lie algebra 
\[ \g = \hcsp(V,\omega) := \heis(V,\omega) \rtimes \csp(V,\omega), \] 
where 
$(V,\omega)$ is a symplectic vector space, 
$\heis(V,\omega) = \R \oplus V$ is the corresponding Heisenberg algebra 
with the bracket $[(z,v),(z',v')] = (\omega(v,v'),0)$, and 
\[ \csp(V,\omega) := \sp(V,\omega) \oplus \R \id_V \] 
is the {\it conformal symplectic Lie algebra} of $(V,\omega)$. 
The hyperplane ideal $\fj := \heis(V,\omega) \rtimes \sp(V,\omega)$ 
(the {\it Jacobi algebra})  
can be identified by the linear isomorphism 
\[ \phi \colon \fj \to \Pol_{\leq 2}(V), \qquad 
\phi(z,v,x)(\xi) := z + \omega(v,\xi) + \frac{1}{2} \omega(x\xi,\xi), \quad 
\xi \in V \] 
with the Lie algebra of polynomials 
$\Pol_{\leq 2}(V)$ of degree $\leq 2$ on $V$, 
endowed with the Poisson bracket (\cite[Prop.~A.IV.15]{Ne99}).  
The set 
\[ C := \{ f \in \Pol_{\leq 2}(V) \colon f \geq 0 \}  \] 
is a pointed generating invariant cone in~$\fj$. 
The element $h_0 := \id_V$ defines a derivation on $\fj$ by 
$(\ad h_0)(z,v,x) = (2z,v,0)$ for $z \in \R, v \in V, x \in \sp(V,\omega)$. 
Any involution $\tau_V$ on $V$ satisfying $\tau_V^*\omega = - \omega$ 
defines by 
\begin{equation}
  \label{eq:invtau}
\tilde\tau_V(z,v,x) := (-z,\tau_V(v), \tau_V x \tau_V) 
\end{equation}
an involution on $\g$ with $\tilde\tau_V(h_0) = h_0$,  
and $-\tilde\tau_V(C) = C$ follows from 
\[ \phi(\tilde \tau_V(z,v,x)) = - \phi(z,v,x) \circ \tau_V.\] 
Considering $\tau_V$ as an element of $\sp(V,\omega)$, 
the element $h := \shalf(\id_V + \tau_V)\in \csp(V,\omega)$ defines 
a $3$-grading of $\g$ because $\ad h$ is diagonalizable 
with eigenvalues $\pm 1,0$. Writing $V = V_1 \oplus V_{-1}$ for 
the $\tau_V$-eigenspace decomposition, we have 
\[ \g^{-1} = 0 \oplus 0 \oplus \sp(V,\omega)^{-1}, \quad 
\g^0 = 0 \oplus V_{-1}\oplus \sp(V,\omega)^0 
\cong V_{-1} \rtimes \gl(V_{-1}), \quad 
\g^1 = \R \oplus V_1\oplus \sp(V,\omega)^1.\] 
Note that 
\[ e^{\pi i \ad h} = (-\tau_V)\,\tilde{}.\] 
Here $\g^1$ can be identified with the space 
$\Pol_{\leq 2}(V_{-1})$ of polynomials of 
degree $\leq 2$ on $V_{-1}$ and 
\[ C_+ = C \cap \g^1 = \{ f \in \Pol_{\leq 2}(V_{-1}) \colon f \geq 0\}.\] 
This  cone is invariant 
under the natural action of the 
affine group $G^0 \cong \Aff(V_{-1})_0 \cong V_{-1} \rtimes \GL(V_{-1})_0$ 
whose Lie algebra is $\g^0$. We also note that 
\[ \g^{-1} = \sp(V,\omega)^{-1} \cong \Pol_2(V_1) \quad \mbox{ and } \quad 
C_- = - C \cap \g^{-1} = \{ f \in \Pol_2(V_1) \colon f \leq 0\}.\] 

Now we turn to the corresponding group and one of its 
irreducible unitary representations. 
Choosing a symplectic basis, we obtain 
an isomorphism with $V \cong V_{-1} \oplus V_1 \cong 
\R^n \oplus \R^n$ with the canonical symplectic 
form specified by 
$\omega((q,0), (0,p)) = \la q,p\ra$ and $\tau_V(q,p) = (-q,p)$. Let 
$\Mp_{2n}(\R)$ denote the {\it metaplectic group}, 
which is the unique non-trivial double cover of $\Sp_{2n}(\R$). 
We consider the group 
\[ G := \Heis(\R^{2n}) \rtimes_\alpha (\R^\times_+ \times \Mp_{2n}(\R)),\] 
where $\R^\times$ acts on $\Heis(\R^{2n}) = \R \times \R^{2n}$ by 
$\alpha_r(z,v) = (r^2 z, rv)$.
Its Lie algebra is $\g = \hcsp(V,\omega)$. 
Then 
\[ \cH 
:= L^2\Big(\R^\times_+, \frac{d\lambda}{\lambda}; L^2(\R^n)\Big) \cong 
= L^2\Big(\R^\times_+ \times \R^n, \frac{d\lambda}{\lambda} \otimes dx\Big),\] 
carries an irreducible representation of $G$, where 
$L^2(\R^n) \cong L^2(V_{-1})$ carries the oscillator representation $U_0$ 
of $\Heis(\R^{2n}) \rtimes \Mp_{2n}(\R)$. 
The Heisenberg group $\Heis(\R^{2n})$ is represented on $\cH$ by 
\begin{align}
(U(z,0,0) f)(\lambda, x) &= e^{i\lambda^2 z} f(\lambda, x), \\ 
(U(0,q,0) f)(\lambda, x) &= e^{i \lambda \la q,x\ra} f(\lambda, x), \\ 
(U(0,0,p) f)(\lambda, x) &= f(\lambda, x - \lambda p).
\end{align}
The group $\Mp_{2n}(\R)$ acts by the metaplectic representation 
on $L^2(\R^n)$ via 
\[ (U(g)f)(\lambda, \cdot ) := U_0(g)f(\lambda, \cdot), \] 
independently of $\lambda$. 
The one-parameter group $\R^\times_+ = \exp(\R h_0)$ acts by 
\[ (U'(r) f)(\lambda,x) := f(r \lambda,x) \quad \mbox{ for }\quad r > 0.\] 
We also note that we have a conjugation $J$ on $\cH$ defined by 
\[ (Jf)(\lambda,x) := \oline{f(\lambda,-x)} \quad \mbox{ satisfying } \quad 
J U(g) J = U (\tau_G(g)),\] 
where $\tau_G$ induces on $\g$ the involution $e^{\pi i \ad h} = (-\tau_V)\,\tilde{}$. 

The positive cone $C_U \subeq \fj$ is the same as the one of the metaplectic 
representation. It intersects $\sp(V,\omega)$ in its unique invariant cone 
of non-negative polynomials of degree $2$ on~$V$. This implies that 
$(C_U)_- = C_-$. To determine $(C_U)_+ = C_U \cap \g^1$, 
we observe that $\g^1$ acts on $L^2(\R^n)\cong L^2(V_-)$  
by multiplication operators. This shows that we also have $(C_U)_+ = C_+$, 
so that we can determine the semigroup $S_\sV$ for the standard subspace 
$\sV \subeq \cH$ with $\Delta_\sV = e^{2\pi i \partial U(h)}$ and $J_\sV = J$. 
It takes the form 
\[ S_\sV = \exp(C_+) G_\sV \exp(C_-),\] 
where $G_\sV = G^0$ is a double cover of $\Aff(\R^n)_0$, its inverse image 
in $\Mp_{2n}(\R)$. 
\end{ex}

\section{Perspectives} 

For an  antiunitary representation $(U,\cH)$ of 
the Lie group $G \rtimes \{\1,\tau_G\}$, 
any element $h \in \g^\tau$ specifies 
a standard subspaces of $\cH$ by the relations 
\[ J_\sV = U(\tau_G) \quad \mbox{ and } \quad 
\Delta_\sV = e^{2\pi i \partial U(h)}.\] 

\subsection{The spaces $\cO_\sV$ and $\cO_h$} 

As we mentioned already in the introduction, 
the $G$-orbit $\cO_\sV := U(G)\sV\cong G/G_\sV$ is a homogeneous space 
on which the inclusion order is invariant, and the order is encoded 
in the semigroup $S_\sV$ by 
\[ U(g_1) \sV \subeq U(g_2) \sV \quad \Leftrightarrow \quad 
g_2^{-1} g_1 \in S_\sV.\] 
The semigroup $S(C_U,h)$ likewise encodes the order on the 
homogeneous space $(\cO_h, \leq_{C_U})$ and the 
Monotonicity Theorem (Theorem~\ref{thm:inc}) asserts that 
$S_\sV \subeq S(h,C_U)$, so that the natural map 
$\pi \colon \cO_\sV \to \cO_h$ is monotone. 
If $\g$ is $3$-graded by $\ad h$ and $\tau = e^{\pi i \ad h}$, then 
$G_\sV$ is an open subgroup of $G_h$ and $\pi$ is a covering 
with $S(h,C_U) = G_h S_\sV$, containing $S_\sV$ as an open subsemigroup, 
so that the concrete ordered 
space $\cO_h$ is a very good model for $(\cO_\sV, \subeq)$. 

\begin{ex}
In general, the connection between $\cO_h$ and $\cO_\sV$ is much less 
intimate, as the example of the $2$-dimensional non-abelian 
Lie group shows. Consider $\g = \R h \oplus \R x$ with 
$[h,x] = \lambda x$ for some $\lambda > 0$. 
Then $C = \R_+ x$ is an invariant cone in $\g$ and the adjoint orbit of $h$ 
is the affine line 
\[ \cO_h = e^{\R \ad x} h = h  + \R x,\] 
endowed with its natural order~$\leq_C$ and 
\[ S(h,C) = \exp(\R h) \exp(\R_+ \lambda x).\] 
If $\lambda \not= 1$, then $C_\pm = \{0\}$ leads 
for representations with $C = C_U$ to $S_\sV = G_\sV$, so that 
the order on $\cO_\sV$ is trivial. Only for $\lambda = \pm 1$ we have 
$S_\sV = S(h,C)$. 
This follows from our result above, but it also can be derived 
directly from the Borchers--Wiesbrock Theorem 
(\cite[\S 3.4]{NO17}).  
\end{ex}

\begin{prob} For a pointed closed convex invariant cone $C \subeq \g$ 
and $h \in \g$, determine the tangent wedge $\L(S(h,C))$ of the semigroup 
$S(h,C)$ in concrete terms. 

Clearly, $h - e^{t\ad x}h \in C$ for $t \geq 0$ implies $[h,x] \in C$, so that 
\[ \L(S(h,C)) \subeq (\ad h)^{-1}(C).\] 
We have seen above that, in the $3$-graded case we have equality because 
\[ (\ad h)^{-1}(C) = C_+ \oplus \g^0 \oplus C_-.\] 
The most important case is when $\ad h$ is real diagonalizable, so that 
$\g = \oplus_{\lambda} \g^\lambda(h)$. 
Then $T_h(\cO_h) \cong [\g, h] = \sum_{\lambda\not=0} \g^\lambda(h)$ 
and $\g^0(h) = \g_h \subeq \L(S(h,C))$. In general it seem rather complicated 
to determine 
\[ (\ad h)^{-1}(C) \cap [h,\g]
= \Big\{ x = \sum_{\lambda \not=0} x_\lambda \colon 
[h,x] = \sum_{\lambda \not=0} \lambda x_\lambda \in C\Big\}.\] 
Only the maximal and minimal 
eigenvalues $\lambda_{\rm min}$ and $\lambda_{\rm max}$ have the property 
that $x  \in C$ implies 
$x_{\lambda_{\rm max}} = \lim_{t \to \infty} e^{-t\lambda_{\rm max}} e^{t \ad h} x \in C,$ 
and likewise $x_{\lambda_{\rm min}} \in C$. 
\end{prob}

\subsection{Covariant nets of standard subspaces} 

As we have seen in Example~\ref{ex:poincare}, 
for the Poincar\'e group $G = P(d)$, 
\[ S(h,C) = \{ g \in G \colon gW_R \subeq W_R\},\] 
so that the ordered space $(\cO_h,\leq_C)$ is isomorphic to the wedge 
space $\cW = G.W_R$ of wedge domains in $\R^{1,d-1}$. 
As such, it provides a natural index set whose elements may be interpreted
as ``special space-time domains''. 

If $\tau$ does not coincide with $e^{\pi i \ad h}$, it is more natural 
to consider pairs $(h,\tau) \in \g \times \Aut(\g)$, where 
$\tau$ is an involution fixing $h$ and to consider $G$-orbits 
$\cO_{(h,\tau)} \subeq \g \times \Aut(\g)$ 
of such pairs. For more on the rich geometric structures 
of such pairs as dilation spaces, we refer to \cite{Ne18}. 

Any pointed convex invariant cone 
$C \subeq \g$ now specifies a natural order on the homogeneous space 
$\cO_{(h,\tau)}$ corresponding to the semigroup 
\[ S = \exp(C_+) G_{(h,\tau)} \exp(C_-), \qquad 
C_\pm := \pm C \cap \g^{-\tau}  \cap \ker(\ad h \mp 1).\] 
Considering the pairs $(h,\tau)$ as abstractions of wedge domains in spacetimes, 
it is now natural to try to classify $G$-covariant maps 
$\cO_{(h,\tau)} \to \Stand(\cH)$ and to study 
the Bisognano--Wichmann property, and their 
causality and duality properties. This project is pursued in \cite{MN20}.  

\subsection{Standard subspaces in Hilbert spaces of distributions}

From the perspective of Quantum Field Theory, it is also interesting 
to see how standard subspaces arise as concrete subspaces of 
Hilbert spaces of distributions. Here one considers a smooth manifold 
$M$ and a positive definite distribution $D$ on $M \times M$, so that 
\[ \la \xi, \eta \ra_D := D(\xi \otimes \oline\eta) \] 
defines a positive semidefinite form on the space 
$C^\infty_c(M,\C)$ of test functions on $M$, hence a Hilbert space 
of distributions $\cH_D \subeq C^{-\infty}(M)$ (cf.~\cite[Ex.~2.4.4]{NO18}). 
For every open subset $\Omega \subeq M$, 
we thus obtain a closed real subspace 
$\sV(\Omega)$ as the closure of the image of $C^\infty_c(\Omega,\R)$. 

We also assume that 
$\alpha \colon \R^\times \to \Diff(M)$ defines an action, such that 
$\alpha(\R^\times_+)$ leaves $D$ invariant and that the involution 
$\tau_M := \alpha(-1)$ satisfies 
\[ \la (\tau_M)_*\xi, (\tau_M)_*\eta \ra_D = \la \eta, \xi \ra_D.\] 
Then we obtain an antiunitary representation $U$ of $\R^\times$, and this 
specifies a standard subspace $\sV \subeq \cH_D$ by 
\[   U(e^t) = \Delta_\sV^{-it/2\pi} \quad \mbox{ for } \quad t \in \R 
\quad \mbox{ and } \quad J_\sV = U(-1).\] 

\begin{problem} Find necessary and sufficient conditions on 
pairs $(\alpha,\Omega)$ such that $\sV = \sV(\Omega)$. 
\end{problem}

This question is studied in \cite{NO20} for the case where 
$M = G$ is a Lie group, $D$ is left invariant 
(hence defined by a positive definite distribution on $G$), and 
$\alpha(e^t)(g) = \exp(th) g \exp(-th)$ for $t \in \R$. 
In this case the semigroups constructed in this article provide natural 
domains on which the real test functions generate a standard subspace. 

In \cite{NOO20} we study the same problem for groups of the form 
$G = (E,+)\rtimes \R^\times$, 
where the Hilbert space $\cH_D$ consists of boundary values 
of holomorphic functions on a tube domain. If $E$ is Minkowski space, then 
our findings show that wedge domains $\Omega\subeq E$ and the corresponding boosts 
provide pairs $(\alpha,\Omega)$ with $\sV_\alpha = \sV(\Omega)$. 

\section*{Acknowledgments}

We thank Gandalf Lechner for an invitation to the 
Simons Center workshop ``Operator Algebras and Applications'' in June 2019, 
where some of the results of this paper have been obtained. 
In particular, we are most grateful 
to Roberto Longo and Gandalf Lechner for pointing 
out that a proof of Theorem~\ref{thm:b.2} 
(Borchers--Wiesbrock Monotonicity) is essentially contained in \cite{Bo00}.
We also thank Konrad Schm\"udgen for illuminating discussions on the subtleties 
of the order on the space of unbounded selfadjoint operators. 

Last, but not least, we also thank Daniel Oeh 
for reading earlier versions of this manuscript.

\appendix

\section{Logarithms of positive operators} 
\mlabel{app:a} 

In  this appendix we collect some background on the order
on the space of not necessarily semibounded selfadjoint operators because 
it is needed in the proof of Theorem~\ref{thm:b.2}.

\begin{defn} (Quadratic form defined by a selfadjoint operator $A$) 
Let $P_A$ denote the spectral measure of $A$ and, for $\xi \in \cH$, write 
$P_A^\xi := \la \xi, P(\cdot) \xi\ra$. Then we define 
\[\cD[A]  
:= \Big\{ \xi \in \cH \colon \int_\R |x|\, dP_A^\xi(x) < \infty \Big\} 
= \cD(|A|^{1/2})\] 
and 
\[ q_A(\xi,\psi) := \int_\R x\, \la \xi, dP_A(x)\psi\ra 
\quad \mbox{ for } \quad \xi, \eta \in \cD[A] \] 
(cf.~\cite[\S VIII.6]{RS80}, \cite[\S 10.2]{Sch12}).  
Clearly, $\cD(A) = \cD(|A|) \subeq \cD[A] = \cD(|A|^{1/2})$, 
but if $A$ is unbounded, then this inclusion is strict. 
\end{defn}

\begin{defn} \mlabel{def:a.2}
For two selfadjoint operators 
$A, B$, semibounded from below, we define 
$A \leq B$ if 
\[ \cD[B] \subeq \cD[A]\quad \mbox{ and } \quad 
q_A(\xi,\xi) \leq q_B(\xi,\xi)\quad \mbox{ for } \quad \xi \in \cD[B]\] 
(\cite[Def.~10.5]{Sch12}). 
If $A$ and $B$ are not semibounded from below, we write $A \preceq B$ if 
\[q_A(\xi,\xi) \leq q_B(\xi,\xi)\quad \mbox{ for } \quad \xi \in 
\cD[A] \cap \cD[B].\] 
\end{defn}

\begin{lem} \mlabel{lem:a.4} For $\Re z > 0$, we have 
\[ \log z = \int_0^\infty  \frac{1}{x+1} - \frac{1}{x+z}\, dx.\] 
\end{lem}

\begin{prf} For $\Re z > 0$, let 
$\gamma(z)(x) := \frac{1}{x+z}$, as a function on 
the half line $(0,\infty)$. Then 
\[ \gamma(z)(x) - \gamma(w)(x) = \frac{w-z}{(x+z)(x + w)} \] 
is integrable over $(0,\infty)$, so that 
$F(z) := \int_0^\infty  \frac{1}{x+1} - \frac{1}{x+z}\, dx$ 
is defined. Next we observe that, for $\Re z > 0$ and 
$|h| < \Re z$, we have 
\[ \frac{\gamma(z+h)(x) - \gamma(z)(x)}{h} + \frac{1}{(x+z)^2} 
= \frac{1}{(x+z)^2} -\frac{1}{(x+z+h)(x + z)} 
= \frac{h}{(x+z)^2(x+z+h)}.\] 
It is easy to see that this expression tends to $0$ in 
$L^1(0,\infty)$ for $h \to 0$. This implies that $F$ is 
holomorphic with 
\[ F'(z) =   \int_0^\infty  \frac{1}{(x+z)^2} \, dx 
=  -\Big[ \frac{1}{x + z}\Big|_0^\infty =  \frac{1}{z}.\] 
As $F(1) = 0$, it follows that $F(z) = \log z$. 
\end{prf}

We want to use the preceding lemma to see that, for a 
selfadjoint operator $A > 0$ ($A \geq 0$ with $\ker A = 0$), we have  
\[ \log(A) = \int_0^\infty  (x + \1)^{-1} - (x + A)^{-1} \, dx\] 
in a suitable sense and derive a suitable version of the operator-monotonicity 
of $\log$ from this 
integral representation. Note that the integrand defines a norm-continuous 
function $(0,\infty) \to B(\cH)$ with a possible singularity in~$0$.

\begin{thm} \mlabel{thm:logmon} If $0 \leq A \leq B$ and $\ker A = 0$, then 
$\ker B = 0$ and 
\[ \log(A) \preceq \log(B), \] 
i.e., $q_{\log(A)}(\xi,\xi) \leq q_{\log(B)}(\xi,\xi)$ for 
$\xi \in \cD[\log(A)] \cap \cD[\log(B)]$.
\end{thm}

\begin{prf}
Let $P_A$ denote the spectral measure of $A$, so  that 
$A = \int_0^\infty x\, dP_A(x).$ 
The condition $\ker A = \{0\}$ means that $P_A(\{0\}) = 0$, so that 
the integral representing $A$ actually extends over the open interval $(0,\infty)$. 
Recall that 
\[ \cD(\log A) = \Big\{ \xi \in \cH \colon \int_0^\infty |\log(x)|^2\, 
dP_A^\xi(x) < \infty \Big\} \] 
and 
\[ \cD[\log A] = \Big\{ \xi \in \cH \colon \int_0^\infty |\log(x)|\, 
dP_A^\xi(x) < \infty \Big\} 
= \cD[\log A_{\geq 1}] \cap \cD[\log A_{< 1}], \] 
where $A_{< 1} = P_A((0,1))A$ and $A_{\geq 1} = P_A([1,\infty)) A$, so that 
$A = A_{< 1} \oplus A_{\geq 1}$. 
We write $\xi \in \cD[\log A]$ accordingly as 
$\xi = \xi_1 \oplus \xi_2$ with $\xi_1 = P_A((0,1))\xi$ and 
$\xi_2 = P_A([1,\infty))\xi$. Then 
we obtain with the Fubini--Tonelli Theorem on iterated integrals 
and Lemma~\ref{lem:a.4}
\begin{align} \label{eq:qtag}
q_{\log(A)}(\xi,\xi) 
&= \int_0^\infty  \log(x)\, dP_A^{\xi}(x) \notag \\
&= \int_{(0,1)} \log(x)\, dP_A^{\xi_1}(x) 
+   \int_{[1,\infty)} \log(x)\, dP_A^{\xi_2}(x)  \notag \\ 
&= \int_{(0,1)} \int_0^\infty \frac{1}{t + 1} - \frac{1}{t + x}\, dt
\, dP_A^{\xi_1}(x) 
+  \int_{[1,\infty)} \int_0^\infty \frac{1}{t + 1} - \frac{1}{t + x}\, dt
\, dP_A^{\xi_2}(x) \notag \\
&=  \int_0^\infty \int_{(0,1)} \frac{1}{t + 1} - \frac{1}{t + x}
\, dP_A^{\xi_1}(x) \, dt
+  \int_0^\infty \int_{[1,\infty)}  \frac{1}{t + 1} - \frac{1}{t + x}
\, dP_A^{\xi_2}(x) \, dt \notag  \\
&=  \int_0^\infty \int_0^\infty  \frac{1}{t + 1} - \frac{1}{t + x}
\, dP_A^{\xi}(x) \, dt \notag \\
&=  \int_0^\infty \la \xi, \big((t + \1)^{-1} - 
(t + A)^{-1}\big) \xi \ra \, dt. 
\end{align}
Here the existence of the latter integral is a consequence of the 
Fubini--Tonelli Theorem. In this sense we have 
\begin{equation}
  \label{eq:formint}
 q_{\log(A)}(\xi,\xi) 
=   \int_0^\infty \big((t + 1)^{-1}\|\xi\|^2 -  q_{(t + A)^{-1}}(\xi,\xi)\big)\, dt 
\quad \mbox{ for }\quad 
\xi \in \cD[\log(A)].
\end{equation}
For $0 < A \leq B$, we have 
\begin{equation}
  \label{eq:14.5}
 -(x + A)^{-1} \leq -(x + B)^{-1} 
\end{equation}
by \cite[Cor.~10.12]{Sch12}, so that \eqref{eq:formint} immediately 
implies the theorem. 
\end{prf}

\begin{cor} If $0 < A \leq B$ are selfadjoint operators, then 
\[ \la \xi, \log(A) \xi \ra \leq \la \xi, \log(B) \xi \ra 
\quad \mbox{ for } \quad \xi \in \cD(\log(A)) \cap \cD(\log(B)).\] 
\end{cor}

\begin{prf} We only have to observe that 
$\cD(\log(A)) \subeq \cD[\log(A)]$ and then use 
Theorem~\ref{thm:logmon}. 
\end{prf}

\begin{rem} Suppose that $c > 0$ and that $A$ and $B$ are selfadjoint 
with $B \geq A \geq c\1$. 
Then $\xi \in \cD[\log(A)]$ 
is equivalent to 
\[ \int_0^\infty\, \log(x)\, dP_A^\xi(x) 
=  \int_c^1 \, \log(x)\, dP_A^\xi(x) + \int_1^\infty \, \log(x)\, dP_A^\xi(x)
< \infty.\] 
In particular, as an element of $\R \cup \{\infty\}$, the integral 
$\int_0^\infty\, \log(x)\, dP_A^\xi(x)$ 
is defined for every $\xi \in \cH$, and in this sense \eqref{eq:qtag} 
and \eqref{eq:14.5} yield  
\begin{align*}
 \int_0^\infty\, \log(x)\, dP_A^\xi(x) 
&=  \int_0^\infty \la \xi, \big((t + \1)^{-1} - 
(t + A)^{-1}\big) \xi \ra \, dt \\
&\leq  \int_0^\infty \la \xi, \big((t + \1)^{-1} - 
(t + B)^{-1}\big) \xi \ra \, dt 
=  \int_0^\infty\, \log(x)\, dP_B^\xi(x)  
\end{align*}
as an equality in $\R \cup \{\infty\}$. 
We conclude that $\cD[\log(B)]\subeq \cD[\log(A)]$, 
so that we also recover from Theorem~\ref{thm:logmon} the well-known 
operator-monotonicity assertion $\log(A) \leq \log(B)$. 

Likewise $\log(B) = - \log(B^{-1})$ shows that 
$0 < A \leq B \leq C \1$ for some $C > 0$ implies that 
$C^{-1} \1 \leq B^{-1} \leq A^{-1}$, so that 
\[ \cD[\log(A)] = \cD[\log(A^{-1})] \subeq 
\cD[\log(B^{-1})] = \cD[\log(B)]\] 
and therefore 
\[ \cD[\log(B)]\cap \cD[\log(A)] = \cD[\log(A)].\] 
This shows that, if $0 \in \Spec(A)$, i.e., $\log(A)$ is not bounded from 
below, then $\log(A) \preceq \log(B)$ is not equivalent to $\log(A) \leq 
\log(B)$ in the sense of Definition~\ref{def:a.2}, but we still have 
$-\log(B) \leq - \log(A)$. 
\end{rem}

\section{Root decomposition} 
\mlabel{app:c} 

In this appendix we recall a few concepts related to root decompositions 
of a finite dimensional Lie algebra $\g$ with respect to a compactly 
embedded Cartan subalgebra. This is used in the proofs 
of Theorems~\ref{thm:1.7} and \ref{thm:1.29}.

Let $\g$ be a finite dimensional Lie algebra and 
$\ft \subeq \g$ be a compactly embedded Cartan subalgebra, 
i.e., the closure of $e^{\ad \ft} \subeq \Aut(\g)$ is compact 
and $\ft$ coincides with its own centralizer: $\ft  = \fz_\g(\ft)$. 
Then we have the root decomposition 
\[ \g_\C = \ft_\C \oplus \bigoplus_{\alpha \in \Delta} 
\g_\C^\alpha, \quad \mbox{ where } \quad 
 \g_\C^\alpha := \{ x \in \g_\C \colon (\forall h \in \ft_\C)\ [h,x]= \alpha(h)
x\} \] 
and 
\[ \alpha(\ft) \subeq i \R \quad \mbox{ for every  root } 
\quad \alpha \in \Delta :=  
\{ \alpha \in \ft_\C^* \setminus \{0\} \colon 
\g_\C^\alpha\not= \{0\}\}.\]
For $x + iy \in \g_\C$ we put $(x+ iy)^* := -x + iy$, so that 
$\g = \{ x \in \g_\C \colon x^* = -x\}$. 
We then have $x_\alpha^* \in \g_\C^{-\alpha}$ for 
$x_\alpha \in \g_\C^\alpha$. 
We call a root $\alpha \in \Delta$ 
\begin{itemize}
\item {\it compact}, if 
there exists an $x_\alpha \in \g_\C^\alpha$ with 
$\alpha([x_\alpha, x_\alpha^*]) > 0$. 
\item {\it non-compact}, if 
there exists a non-zero $x_\alpha \in \g_\C^\alpha$ with 
$\alpha([x_\alpha, x_\alpha^*]) \leq 0$. 
\item {\it non-compact simple}, if 
there exists a non-zero $x_\alpha \in \g_\C^\alpha$ with 
$\alpha([x_\alpha, x_\alpha^*]) < 0$. 
\end{itemize}
We write 
$\Delta_k, \Delta_p, \Delta_{p,s} \subeq \Delta$ for the subset of compact, 
non-compact, resp., non-compact simple roots. 
A subset $\Delta^+\subeq \Delta$ is called a {\it positive system} 
if there exists an $x_0 \in \ft$ with $\alpha(x_0) \not=0$ for every 
$\alpha \in \Delta$ and 
\[ \Delta^+  = \{ \alpha \in \Delta \colon i\alpha(x_0) > 0 \}.\] 
A positive system $\Delta^+$ is said to be {\it adapted} 
if $i\alpha(x_0) > i \beta(x_0)$ for 
$\alpha \in \Delta_p^+$ and $\beta \in \Delta_k$ 
(cf.~\cite[Prop.~VII.2.12]{Ne99}). 
To an adapted positive system $\Delta^+$, we associate the cone 
\begin{equation}
  \label{eq:maxcon}
C_{\rm max} := C_{\rm max}(\Delta_p^+) := (i\Delta_p^+)^\star.
\end{equation} 
Now $W_{\rm max} := \oline{\Ad(G)C_{\rm max}}$ is a closed convex invariant cone 
with $W_{\rm max}^0 = \Ad(G) C_{\rm max}^0$. We also note that 
$W_{\rm max} \cap \ft = C_{\rm max}$ (\cite[Lemma~VIII.3.22, 27]{Ne99}).

\section{Projections onto graphs} 
\mlabel{app:b}

For the sake of completeness, we include here some arguments from 
\cite[\S II.1]{Bo00} that are used in the proof of the Monotonicity Theorem 
(Theorem~\ref{thm:b.2}). 

\begin{lem}    \mlabel{lem:b.1} 
Let $\cH_1$ and $\cH_2$ be complex Hilbert spaces, 
let $S \colon \cH_1 \supeq \cD(S) \to \cH_2$ be a closed operator from $\cH_1$ 
to $\cH_2$, and let $P = \pmat{p_{11} & p_{12} \\ p_{21} & p_{22}}$ 
denote the orthogonal projection onto the closed subspace 
$\Gamma(S) = \{ (x,Sx) \in \cH_1 \oplus \cH_2\colon x \in \cD(S)\}$, 
written as a $(2 \times 2)$-matrix. 
Then 
\[ p_{11} = (\1 + S^*S)^{-1} \quad \mbox{ and } \quad 
p_{12}\res_{\cD(S^*)} = (\1 + S^*S)^{-1}S^*.\] 
\end{lem}

\begin{prf}
The relation $P = P^2$ implies $p_{11} = p_{12}p_{21} + p_{11}^2$. 
As $P(\xi,0) = (p_{11}\xi, S p_{11} \xi) \in \Gamma(S)$, we have 
\[ (\xi,0) - (p_{11}\xi, S p_{11} \xi) 
= ((\1 - p_{11})\xi, - S p_{11}\xi) \bot \Gamma(S).\] 
Further, 
\begin{equation}
  \label{eq:graphorth}
 \Gamma(S)^\bot = \{ (-S^* \psi, \psi) \colon \psi \in \cD(S^*) \}
\end{equation}
now shows that 
$\1 - p_{11} = S^*S p_{11}$, i.e., 
$\1= (\1 + S^*S) p_{11}.$ 
As $\1 + S^*S$ is injective, it follows that 
$p_{11} = (\1 + S^*S)^{-1}$. 

From $P(0,\xi) = (p_{12}\xi, S p_{12} \xi) \in \Gamma(S)$, we likewise get 
\[ (0,\xi) - (p_{12}\xi, S p_{12} \xi) 
= (- p_{12}\xi, (\1 - S p_{12})\xi) \bot \Gamma(S).\] 
With \eqref{eq:graphorth}, this leads to 
$p_{12} = S^* (\1 - S p_{12})$. For $\xi \in \cD(S^*)$, we thus obtain 
$p_{12}\xi = S^*\xi - S^*S p_{12}\xi$, so that 
$(\1 + S^*S) p_{12}\xi = S^* \xi$, and finally 
$p_{12} \xi = (\1 + S^*S)^{-1} S^*\xi$. 
\end{prf}

Lemma~\ref{lem:b.1} can be used to characterize operators 
on $\cD(S)$ which are bounded in the graph topology. 

\begin{lem} \mlabel{lem:b.2} Let $\cH_3$ be a Hilbert space and 
$A \colon \cD(S) \to \cH_3$ be a linear map. 
Then $A$ is continuous with respect to  the graph topology on $\cD(S)$ 
if and only if the operators 
\[ A (\1 + S^*S)^{-1} \colon \cH_1 \to \cH_3 
\quad \mbox{ and } \quad A (\1 + S^*S)^{-1} S^* \colon \cD(S^*) \to \cH_3\] 
are bounded, where $\cD(S^*)\subeq \cH_2$ carries the subspace topology.   
\end{lem}

\begin{prf} The operator 
$A$ is continuous in the graph topology if and only if the operator 
\[ \tilde A \colon \Gamma(S) \to \cH_3, \quad (\xi, S \xi) \mapsto A \xi \]
is bounded, and this is equivalent to the boundedness of 
$\tilde A \circ P  \colon \cH_1 \oplus \cH_2 \to \cH_3.$ 
As this operator has the form
$\tilde A P(\xi_1, \xi_2) = A(p_{11}\xi_1 + p_{12}\xi_2),$ 
its continuity is by Lemma~\ref{lem:b.1} equivalent to the boundedness of 
\[ Ap_{11} = A (\1 + S^*S)^{-1} \quad \mbox{ and } \quad 
A p_{12}\res_{\cD(S^*)} = A (\1 + S^*S)^{-1} S^*.\qedhere\] 
\end{prf}

\begin{rem}
If $S \colon \cH \supeq \cD(S) \to \cH$ is an antilinear operator 
and $\cH^{\rm op}$ denotes $\cH$, endowed with the opposite complex structure,
then Lemma~\ref{lem:b.1} 
applies with $\cH_1 = \cH$ and $\cH_2 = \cH^{\rm op}$. 
\end{rem}

\end{document}